\documentclass[a4paper,11pt]{article}
\usepackage[T1]{fontenc}

\usepackage{amssymb}

\usepackage{amsmath,dsfont}

\usepackage{graphicx}

\usepackage{natbib}

\newtheorem{lem}{Lemma}

\newtheorem{thm}{Theorem}

\newtheorem{coro}{Corollary}

\newcommand{\X}{{\cal X}}

\newcommand{\R}{{\mathbb{R}}}

\newcommand{\Eb}{\text{E}}
\newcommand{\Pb}{\text{P}}

\date{}

\begin{document}
\title{Recursive estimation of nonparametric regression with functional covariate.}
\author{Aboubacar AMIRI$^\ddag$ 
\and Christophe CRAMBES\footnote{Corresponding author: christophe.crambes@univ-montp2.fr}
\footnote{Universit\'e Montpellier 2, Institut de Math\'ematiques et de Mod\'elisation de Montpellier, Place Eug\`ene Bataillon,  34090 Montpellier, France. Email: christophe.crambes@univ-montp2.fr} 
\and Baba THIAM\footnote{Universit\'e Lille Nord de France, Universit\'e Lille 3, Laboratoire EQUIPPE,
EA 4018, Villeneuve d'Ascq, France. Emails: aboubacar.amiri@univ-lille3.fr, baba.thiam@univ-lille3.fr}
}
\maketitle

\begin{abstract}
\noindent The main purpose is to estimate the regression function of a real random variable with functional explanatory variable by using a recursive nonparametric kernel approach. The mean square error and the almost sure convergence of a family of recursive kernel estimates of the regression function are derived. These results are established with rates and precise evaluation of the constant terms. Also, a central limit theorem for this class of estimators is established. The method is evaluated on simulations and real data set studies.
 \end{abstract}

\noindent{\it Keywords~:} Functional data, recursive kernel estimators, regression function, quadratic error, almost sure convergence, asymptotic normality.\\
\noindent{\it MSC~:} 62G05,~~62G07,~~62G08.

\section{Introduction}\label{intro}
Functional data analysis is a branch of statistics that has been studied frequently and developed over recent years. This type of data appears in many practical situations such as continuous phenomena. Thus, the possible application fields that favor the use of functional data are broad and include the following: climatology, economics, linguistics, medicine, and so on. Because the pioneering work of authors such as \cite{RamsayDalzell} and  \citet{FrankFriedman}, many developments have been investigated to build theories and methods regarding functional data; for instance, how can the mean or the variance of functional data be defined? What type of model can be considered with functional data? These papers also emphasize the drawback of merely using multivariate methods with this type of data; rather, they suggest considering these data as objects that belong to a specific functional space. The monographs of \cite{RamsaySilverman05, RamsaySilverman02} provide and overview of both the theoretical and practical aspects of functional data analysis.

Regression is one of the most studied functional models. In this model, the variable of interest $Y$ is real and the covariate $\mathcal{X}$ belongs to a functional space $\mathcal{E}$ endowed with a semi-norm $\left\| \cdot\right\|$. Thus, the regression model is written
\begin{equation}
Y = r(\mathcal{X}) + \varepsilon, \label{regression}
\end{equation}

\noindent where $r: \mathcal{E} \longrightarrow \mathbb{R}$ is an operator and $\varepsilon$ is the random error term. Much research has investigated this model when the operator $r$ is supposed to be linear, which has contributed to the popularity of the so-called functional linear model. In this linear context, the operator $r$ is written as $\langle \alpha , . \rangle$ where $\langle . , . \rangle$ denotes the inner product of the space $\mathcal{E}$ and $\alpha$ belongs to $\mathcal{E}$. Thus, the goal is to estimate the unknown function $\alpha$. We refer the reader to \cite{CFS03} and \cite{CKS09} for different methods of estimating $\alpha$. Another way to approach the (\ref{regression}) model is to think in a nonparametric way. Many authors have also investigated this direction. Recent advances on the topic have been the subject of a bibliographical review in \cite{MantVieu} and monographs by \citet{Ferraty2006}, \citet{FerratyRomain}, thereby providing the theoretical and practical properties of a kernel estimator of the operator $r$. Specifically, if $\left( \mathcal{X}_{i} , Y_{i} \right)_{i=1, \ldots, n}$ is a sample of independent and identically distributed couples with the same law as $(\mathcal{X} , Y)$, then, this kernel estimator is defined for all $\chi \in \mathcal{E}$ by

\begin{equation}
r_{n} (\chi) := \frac{{\displaystyle \sum_{i=1}^{n} Y_{i} K \left( \frac{\left\| \chi - \mathcal{X}_{i} \right\|}{h} \right)}}{{\displaystyle \sum_{i=1}^{n} K \left( \frac{\left\| \chi - \mathcal{X}_{i} \right\|}{h} \right)}}, \label{estimFV}
\end{equation}

\noindent where $K$ is a kernel and $h > 0$ is a bandwidth. \cite{Masry2005} considered the asymptotic normality of \eqref{estimFV} in the dependent case, while \cite{LINGWU} obtained almost sure convergence. This nonparametric regression estimator raises several problems because the choice of the semi-norm $\left\| \cdot \right\|$ of the space $\mathcal{E}$ and the choice of the bandwidth, \ldots. With regard to bandwidth, many solutions have been considered when the covariate is real (e.g., cross validation). Recently, \citet{Amiri} studied an estimator using a sequence of bandwidths in the multivariate setting that allowed the computation of this estimator in a recursive way thereby generalizing previous research (\citet{DW80} and  \citet{AL76}). This estimator has acceptable theoretical properties (from the point of view of the mean square error and almost sure convergence). It is also of practical interest: for instance, it presents a computational time gain when researchers want to predict new values of the variable of interest when new observations appear. This case is not true for the basic kernel estimator, which must be computed again using the whole sample. The purpose of the current study is to adapt the recursive estimator studied in \cite{Amiri} to a case in which the covariate is functional.

The remainder of the paper is organized as follows. Section 2 defines the recursive estimator of the operator $r$ when the covariate $\mathcal{X}$ is functional and presents the asymptotic properties of this estimator. Section 3 evaluates the performance of our estimator using a simulation study and a real dataset. Finally, the proofs of the theoretical results are presented in Section 4.

\section{ Functional regression estimation}
\subsection{Notations and assumptions}
Let $({\cal X},Y)$ be a pair of random variables defined in $\left(\Omega,\mathcal{A},P \right),$   with values on ${\cal E}\times\R ,$ where ${\cal E}$ is a Banach space endowed with a semi-norm $\left\| \cdot\right\|$.
Assume that $({\cal X}_i,Y_i)_{i=1,\ldots,n} $ is a sample of $n$ random variables independent and identically distributed,   having  the same distribution as $({\cal X},Y)$.  The  \eqref{regression} model  is then rewritten as   $$Y_i=r({\cal X}_i)+\varepsilon_i,~~i=1,\ldots,n,$$
where  for any $i=1,\ldots,n$, $\varepsilon_i$ is a  random variable  such that   $\text{E}(\varepsilon_i|{\cal X}_i)=0$ and  $\text{E}(\varepsilon_i^2|{\cal X}_i)=\sigma_\varepsilon^2({\cal X}_i)<\infty.$\\
 Nonparametric regression aims to estimate the functional  $r(\chi):=\text{E}\left( Y |{\cal X}=\chi \right),$ for $\chi \in {\cal E}$. To this end, let us consider the family of recursive estimators indexed by a parameter $\ell\in[0,1],$ and defined by
 $$r_n^{[\ell]}({\chi}) :=\frac{\sum\limits_{i=1}^n\frac{Y_i}{F(h_i)^{\ell}}K\left( \frac{\|{\chi}-{\cal X}_i\|}{h_i}\right)}
{\sum\limits_{i=1}^n\frac{1}{F(h_i)^{\ell}}K\left( \frac{\|{\chi}-{\cal X}_i\|}{h_i}\right)},$$
where $K$ is a kernel, $(h_n)$ a sequence of bandwidths and $F$ the cumulative distribution function of the random variable $\|\chi-{\cal X}\|$.   Our family of estimators  is a recursive modification of the estimate defined in \eqref{estimFV}  and can be computed recursively by
$$r_{n+1}^{[\ell]}(\chi)=\frac{\left[\sum\limits_{i=1}^nF(h_i)^{1-\ell}\right]\varphi_n^{[\ell]}(\chi)+\left[\sum\limits_{i=1}^{n+1}F(h_i)^{1-\ell}\right]Y_{n+1}K_{n+1}^{[\ell]}\left(\|\chi-{\cal X}_{n+1} \|\right)}{\left[\sum\limits_{i=1}^nF(h_i)^{1-\ell}\right]f_n^{[\ell]}(\chi)+\left[\sum\limits_{i=1}^{n+1}F(h_i)^{1-\ell}\right]K_{n+1}^{[\ell]}\left(\|\chi-{\cal X}_{n+1}\| \right)},$$
with

\begin{eqnarray}
\varphi_n^{[\ell]}({\chi})\label{phin}=\frac{\sum\limits_{i=1}^n\frac{Y_i}{F(h_i)^{\ell}}K\left( \frac{\|{\chi}-{\cal X}_i\|}{h_i}\right)}{\sum\limits_{i=1}^nF(h_i)^{1-\ell}}, ~
f_n^{[\ell]}({\chi})\label{fn}=\frac{\sum\limits_{i=1}^n\frac{1}{F(h_i)^{\ell}}K\left( \frac{\|{\chi}-{\cal X}_i\|}{h_i}\right)}{\sum\limits_{i=1}^nF(h_i)^{1-\ell}},
\end{eqnarray}
and $K_i^{[\ell]}(\cdot):=\frac{1}{F(h_i)^{\ell}\sum\limits_{j=1}^iF(h_j)^{1-\ell}}K\left(\frac{\cdot}{h_i}\right)$. \\
More precisely, $r_n^{[\ell]}({\chi})$ is the adaption to the functional model of the finite-dimensional recursive family of estimators introduced by \cite{Amiri}, which includes the well-known recursive ($\ell=0$) and semi recursive ($\ell=1$) estimators. The recursive property of this estimator class is clearly useful in sequential investigations and for large sample sizes because the addition of a new observation means that the non-recursive estimators must be recomputed. In addition, we need to store extensive data to calculate them.\\

We assume that the following assumptions hold:
\begin{description}
\item [H1] The operators  $ r $ and $ \sigma_\varepsilon^2$ are continuous on a neighborhood of ${\chi}$, and $F(0)=0$. Moreover, the function
$\varphi(t):=\text{E}\left[\{r({\cal X})-r(\chi)\}|\|{\cal X}-\chi\|=t\right]$ is assumed to be derivable at $t=0$.
\item [H2] $K$ is nonnegative bounded kernel with support on the compact $[0,1]$ such that $\inf\limits_{t\in [0,1]}K(t)>0$.
\item [H3] For any $s\in[0,1], \tau_h(s):=\frac{F(hs)}{F(h)}\rightarrow\tau_0(s)<\infty$ as $h\rightarrow0$.
\item [H4]
\begin{description}
 \item[(i)] $h_n\rightarrow 0, nF(h_n)\rightarrow\infty$ and $A_{n,\ell}:=\displaystyle\frac{1}{n}\sum\limits_{i=1}^n\frac{h_i}{h_n}\left[\frac{F(h_i)}{F(h_n)}\right]^{1-\ell}\rightarrow\alpha_{[\ell]}<\infty$ as $n \to \infty$.
\item[(ii)] $\forall r\leq2$, $B_{n,r}:=\displaystyle\frac{1}{n}\sum\limits_{i=1}^n\left[\frac{F(h_i)}{F(h_n)}\right]^r\rightarrow\beta_{[r]}<\infty,$ as $n \to \infty$.
\end{description}
 \end{description}
     Assumptions ${\bf H1}$, ${\bf H2}$ and the first part of ${\bf H4}$ are common in nonparametric regression. They have been used by \cite{Ferraty2007} and are the same as those classically used in finite-dimensional settings. The conditions $A_{n,\ell}\rightarrow\alpha_{[\ell]}<\infty$ and ${\bf H4(ii)}$ are particular to the recursive problem and the same as those used in finite-dimensional cases.
Note that $F$ plays a crucial role in our calculus: Its limit at  zero, and a fixed ${\cal \chi}$ is known as `small ball probability'. Before announcing our results, let us provide typical examples of bandwidths and small ball probabilities that satisfy ${\bf H3}$ and  ${\bf H4}$ (see \cite{Ferraty2007} for more details).\\
 If ${\cal X}$ is a fractal (or geometric) process, then the small ball probabilities are of the form
$F(t)\sim c'_{\chi}t^\kappa, $
where  $c'_{\chi}$ and $\kappa$ are positive constants, and $\left\| \cdot\right\|$ might be a supremum norm, an $L^p$ norm or a Besov norm.  The choice of bandwidth $h_n=An^{-\delta}$ with $A>0 ~ \text{and}~0<\delta<1$ implies $F(h_n)=c''_\chi n^{-\delta\kappa }$, $c''_\chi>0$. Thus, ${\bf H3}$ and ${\bf H4}$ hold when $\delta\kappa<1$.  In fact, assumption ${\bf H3}$ and the first part of ${\bf H4}$ are clearly non-restrictive because they are the same as those used in  the non-recursive case. With regard to ${\bf H4 (ii)}$, if  $\delta\kappa r<1$, then $\sum\limits_{i=1}^ni^{-\delta\kappa r }\sim\frac{n^{1-\delta\kappa r}}{1-\delta\kappa r}$. Thus, the condition is satisfied when $\beta_{[r]}=\frac{1}{1-\delta\kappa r}$. The same argument is also valid for $A_{n,\ell}$, if $\max\{\kappa r, 1+\kappa(1-\ell)\}<1/\delta$.

\subsection{Main  results}
Let us introduce the following notations from \cite{Ferraty2007}:
\begin{eqnarray*}
M_0&=&K(1)-\int_{0}^1(sK(s))'\tau_0(s)ds,~M_1=K(1)-\int_{0}^1K'(s)\tau_0(s)ds,\\
M_2&=&K^2(1)-\int_{0}^1(K^2(s))'\tau_0(s)ds.\\
\end{eqnarray*}
 We can establish the asymptotic mean square error of our recursive estimate.

\begin{thm}\label{Bias_Var}
Under the assumptions  ${\bf H1}-{\bf H4}$, we have
\begin{eqnarray*}
E\left[r_n^{[\ell]}({\chi})\right]-r(\chi)&=&\varphi'(0)\frac{\alpha_{[\ell]}}{\beta_{[1-\ell]}}\frac{M_0}{M_1}h_n[1+o(1)]+O\left[\frac{1}{nF(h_n)}\right], \\
Var\left[r_n^{[\ell]}({\chi})\right]&=& \frac{\beta_{[1-2\ell]}}{\beta_{[1-\ell]}^2}\frac{M_2}{M_1^2}\sigma_\varepsilon^2(\chi)\frac{1}{nF(h_n)}\left[1+o(1)\right].
\end{eqnarray*}
\end{thm}
Theorem \ref{Bias_Var} is an extension of \cite{Ferraty2007}  result to the class of recursive estimators.  Using a bias-variance representation and an additional condition, the asymptotic mean square error of our estimators is established in the following result:

\begin{coro}\label{MSE}
Assume that the assumptions of Theorem \ref{Bias_Var} hold. If there exists a constant $c>0$ such that $nF(h_n)h_n^2\to c,$ as $n\rightarrow \infty,$ then
\begin{eqnarray*}
\lim_{n\to \infty}nF(h_n)\Eb\left[\left(r_n^{[\ell]}({\chi})-r(\chi)\right)^2\right]&=&\left[\frac{\beta_{[1-2\ell]}}{\beta_{[1-\ell]}^2}\frac{M_2\sigma_\varepsilon^2(\chi)}{M_1^2}+\frac{c\alpha_{[\ell]}^2}{\beta_{[1-\ell]}^2}\frac{\varphi'(0)^2M_0^2}{M_1^2}\right].
\end{eqnarray*}
In particular, if
$\X$ is fractal (or geometric process) with $F(t)\sim c'_{\chi}t^\kappa,$ then  the choice $h_n=An^{-\frac{1}{\kappa+2}}$,  $A,\kappa>0$, implies that
\begin{eqnarray*}
\lim_{n\to \infty}n^{\frac{2}{2+\kappa}}E\left[\left(r_n^{[\ell]}({\chi})-r(\chi)\right)^2\right]&=&\left[\frac{\beta_{[1-2\ell]}}{\beta_{[1-\ell]}^2}\frac{M_2\sigma_\varepsilon^2(\chi)}{c'_\chi A^\kappa M_1^2}+\frac{\alpha_{[\ell]}^2}{\beta_{[1-\ell]}^2}\frac{\varphi'(0)^2M_0^2A^2}{M_1^2}\right].
\end{eqnarray*}
\end{coro}
 \cite{bosq-cheze} established a similar result for the Nadaraya-Watson estimator using the finite-dimensional setting and continuous time processes.\\

To obtain the almost sure convergence rate of our estimator, we assume that the following additional assumptions hold.
\begin{description}
\item [H5]There exist   $\lambda>0$ and $\mu>0$  such that $\Eb\left[\exp\left(\lambda |Y|^\mu\right)\right]<\infty.$
\item [H6] $
\lim\limits_{n\rightarrow+\infty}
\frac{nF(h_n)(\ln n)^{-1-\frac{2}{\mu}}}{\left(\ln\ln n\right)^{2(\alpha+1)}}=\infty  \text{ for some } \alpha \geq0 \text{ and } \lim\limits_{n\rightarrow+\infty}(\ln n)^{\frac{2}{\mu}}F(h_n)=0$.
\end{description}

 Assumption  ${\bf H5}$ is clearly met if $Y$ is bounded, which implies that
\begin{eqnarray}\label{moment-expo}\text{E}\left( \max\limits_{1\leq i\leq n}|Y_i|^p\right) = O[(\ln n)^{p/\mu}], \forall  p \geq1,n \geq 2.\end{eqnarray}
In fact, if we set $M=\left\{\begin{array}{cc}\left(\frac{p-\mu}{\lambda\mu}\right)^{1/\mu}&\text{ if } p>\mu \\
0 &\text{ else }\end{array} \right.$, then one is able to write:
\begin{eqnarray*}\Eb\left(\max_{1\leq i\leq n}|Y_i|^p\right)\leq M^p+\Eb\left(\max_{1\leq i\leq n}|Y_i|^p{\mathds {1}}_{\{|Y_i|>M\}}\right).\end{eqnarray*}
Because the function  $x\mapsto(\ln x)^{p/\mu}$ for all $p\geq 1$ is concave down on set $] \max\{1,\exp(\frac{p}{\mu} -1)\}, +\infty[$, Jensen's inequality (via assumption  ${\bf H5}$) implies that
\begin{eqnarray*} \Eb\left(\max_{1\leq i\leq n}|Y_i|^p{\mathds{1}}_{\{|Y_i|>M\}}\right)\leq&\left[\ln \Eb\exp\left(\lambda\max\limits_{1\leq i\leq n}|Y_i|^\mu{\mathds{1}}_{\{|Y_i|>M\}}\right)\right]^{p/\mu}\\
  \leq& \left[\ln \sum\limits_{i=1}^n\Eb\exp\left(\lambda|Y_i|^\mu\right)\right]^{p/\mu}=O[(\ln n)^{p/\mu }], \end{eqnarray*}
and (\ref{moment-expo}) follows.  An example of a sequence of random variables, $Y_i$ satisfying ${\bf H5}$ (and then (\ref{moment-expo})) is the standard Gaussian distribution, with $\lambda=1$ and $\mu=2$.  \cite{bosq-cheze} used relation (\ref{moment-expo})  in the multivariate framework to establish the optimal quadratic error of the Nadaraya-Watson estimator.  Assumption ${\bf H6}$ is satisfied when $\X$ is a fractal or non-smooth, whereas condition $\lim\limits_{n\to \infty}F(h_n)(\ln n)^{\frac{2}{\mu}}=0$ is not necessary when $\mu\geq2$.\\
We can write the following theorem for our estimator of the regression operator.
\begin{thm}\label{Cv_ps_rnl}
     Assume that ${\bf H1}-{\bf H6}$ hold.
 If   $ \displaystyle\lim\limits_{n\rightarrow+\infty}nh_n^{2}=0$, then
$$\limsup_{n\rightarrow\infty} \left[\frac{nF(h_n)}{\ln\ln~n}\right]^{1/2}\left[r_n^{[\ell]}(\chi)-r(\chi)\right]=\frac{\left[2\beta_{[1-2\ell]}\sigma^2_\varepsilon(\chi)M_2\right]^{1/2}}{\beta_{[1-\ell]}M_1} ~ a.s.$$
\end{thm}

 The choices of bandwidths and small ball probabilities previously provided are typical examples that satisfy the condition $ \displaystyle\lim\limits_{n\rightarrow+\infty}nh_n^{2}=0.$
Case $\ell=1$ of Theorem \ref{Cv_ps_rnl} is an extension to the functional setting of \cite{Roussas92} result  concerning the almost sure convergence of Devroye-Wagner's estimator.
Note, that in the non-recursive framework, the rate of convergence obtained is of the form $\left[\dfrac{nF(h_n)}{\ln n}\right]^{1/2}$ (see Lemma 6.3 in \cite{Ferraty2006}). Unlike the non-recursive case, the rate of convergence of the recursive estimators are obtained with exact upper bounds.\\
To obtain asymptotic normality, we make the following additional assumption, which was clearly verified by the choices of bandwidths and the small ball probabilities above.
\begin{description}
\item [H7] For any $\delta>0$, $\lim\limits_{n\to \infty}\dfrac{(\ln n)^\delta}{\sqrt{nF(h_n)}}=0$.
\end{description}
\begin{thm}\label{normality}
Assume that ${\bf H1}-{\bf H5}$ and ${\bf H7}$ hold. If there exists $c\geq 0$ such that $\lim\limits_{n \to \infty}h_n\sqrt{nF(h_n)}=c$,
 then
\begin{eqnarray*}
\sqrt{nF(h_n)}\left(r_n^{[\ell]}(\chi)-r(\chi)\right)\stackrel{\mathcal{D}}{\rightarrow}\mathcal{N}\left(c\frac{\alpha_{[\ell]}}{\beta_{[1-\ell]}}\frac{M_0}{M_1}\varphi'(0),\ \ \frac{\beta_{[1-2\ell]}}{\beta_{[1-\ell]}^2}\frac{M_2}{M_1^2}\sigma_\varepsilon^2(\chi)\right).
\end{eqnarray*}
\end{thm}
 Note that the choices of bandwidths and small ball probabilities above imply that  ${\beta_{[1-2\ell]}}/{\beta_{[1-\ell]}^2}<1.$ Thus, the recursive estimators are more efficient than classical estimators in this case in the sense that their asymptotic variance is small for a given $M_1, M_2$ and $\sigma_\epsilon(\chi)$.\\
Because the result of Theorem \ref{normality} depends on unknown quantities, we derived a usable asymptotic distribution for the following corollary of case $\ell=0$.
\begin{coro}\label{normality2}
Assume that ${\bf H1}-{\bf H5}$ and ${\bf H7}$ hold. If $\lim\limits_{n \to \infty}h_n\sqrt{nF(h_n)}=0$ and  $\lim\limits_{n \to \infty}\frac{\ln n}{nF^2(h_n)}=0$, then for any consistent estimators $\widehat{M}_i, i =1,2$ and $\widehat\sigma_\varepsilon(\chi)$ of $M_i, i=1,2$ and $\sigma_\varepsilon(\chi)$ respectively, we have
\begin{eqnarray*}
\sqrt{n\widehat{F}(h_n)}\sqrt{\frac{\widehat\beta_{[1]}\widehat{M_1}^2}{\widehat{M}_2\widehat\sigma_\varepsilon^2(\chi)}}\left(r_n^{[0]}(\chi)-r(\chi)\right)\stackrel{\mathcal{D}}{\rightarrow}\mathcal{N}\left(0,1\right),
\end{eqnarray*}
where $\widehat{F}$ and $\widehat\beta_{[1]}$ are the empirical counterparts of $F$ and $\beta_{[1]}$ defined by 
\begin{eqnarray*}
\widehat{F}(t)=\frac{1}{n}\sum_{i=1}^n\mathds{1}_{\{\|\X_i-\chi\|\leq t\}}\text{ and }
\widehat\beta_{[1]}=\displaystyle\frac{1}{n}\sum\limits_{i=1}^n\frac{\widehat F(h_i)}{\widehat F(h_n)}.
\end{eqnarray*}
\end{coro} 
Corollary \ref{normality2} is similar to the result obtained by \cite{Ferraty2007} in the non-recursive case. The assumptions used to establish this result are fulfilled by the choices of bandwidths and small ball probabilities above. If we consider the uniform kernel $K(\cdot)=\mathds{1}_{[0,1]}(\cdot)$, then the asymptotic distribution in Corollary \ref{normality2} can be rewritten as 
\begin{eqnarray*}
\sqrt{\frac{n\widehat{F}(h_n)\widehat\beta_{[1]}}{\widehat\sigma_\varepsilon^2}}\left(r_n^{[0]}(\chi)-r(\chi)\right)\stackrel{\mathcal{D}}{\rightarrow}\mathcal{N}\left(0,1\right).
\end{eqnarray*}
According to Corollary \ref{normality2}, the asymptotic confidence band of $r(\chi)$ with level $1-\alpha$ is given by
$$ \left[ r_n^{[0]}(\chi) - z_{1-\alpha/2} \sqrt{n\widehat{F}(h_n)}\sqrt{\frac{\widehat\beta_{[1]}\widehat{M_1}^2}{\widehat{M}_2\widehat\sigma_\varepsilon^2(\chi)}} ; r_n^{[0]}(\chi) + z_{1-\alpha/2} \sqrt{n\widehat{F}(h_n)}\sqrt{\frac{\widehat\beta_{[1]}\widehat{M_1}^2}{\widehat{M}_2\widehat\sigma_\varepsilon^2(\chi)}} \right], $$
where $z_{1-\alpha/2}$ is the quantile of order $1-\alpha/2$ of the standard normal distribution.

\section{Simulation study and real dataset example}
To observe the behavior of our recursive estimator in practice, this section considers a simulation study. We simulated our data in the following way: The curves $\mathcal{X}_{1}, \ldots$, and $\mathcal{X}_{n}$ are standard Brownian motions on $[0,1]$, with $n = 100$. Each curve is discretized into $p = 100$ equidistant points on $[0,1]$. The operator $r$ is defined by ${\displaystyle r (\chi) = \int_{0}^{1} \chi (s)^{2} \, ds}$. The error $\varepsilon$ is simulated as a Gaussian random variable with mean $0$ and standard deviation $0.1$. The simulations were repeated $500$ times to compute the prediction errors for a new curve $\chi$, which was also simulated as a standard Brownian motion on $[0,1]$.\\
In this functional context, the estimator depends on the choice of many parameters: the semi-norm $\left\| \cdot \right\|$ of the functional space $\mathcal{E}$, the sequence of bandwidths $(h_{n})$, the kernel $K$, the parameter $\ell$ and the distribution function $F$ in case $\ell\neq0$. Because the choice of kernel $K$ is not crucial, we used the quadratic kernel defined by $K(u) = \left( 1 - u^{2} \right) \mathds{1}_{[0,1]} (u)$ for all $u \in \mathbb{R}$, which behaves correctly in practice and is easy to implement. We estimated the distribution function $F,$ using the empirical distribution function, which is uniformly convergent.
\subsection{Choice of the bandwidth}
In this simulation, the semi-norm was based on the principal components analysis of the curves that retained $3$ principal components (see\cite{BCF97} for a description of this semi-norm), whereas $\ell$ is fixed and equal to $0$. We show below that this parameter $\ell$ has a negligible influence on the behavior of the estimator.\\
 We chose to take a sequence of bandwidths $ h_{i} = C \max\limits_{i = 1, \ldots, n} \left\| \mathcal{X}_{i} - \chi \right\| i^{-\nu}$, $ \text{ for }  1\leq i \leq n, $ with $C\in  \left\{ 0.5, 1, 2, 10 \right\}$ and $\nu\in \left\{ \frac{1}{10} , \frac{1}{8} , \frac{1}{6} , \frac{1}{5} , \frac{1}{4} , \frac{1}{3} , \frac{1}{2} ,  1 \right\}$. \\
At the same time, we also computed the estimator \eqref{estimFV} introduced by \cite{Ferraty2006}.  Following \cite{RachdiVieu}, we introduced an automatic selection of the bandwidth with a cross validation procedure. We used this procedure for the estimator of \cite{Ferraty2006}. For the recursive estimator, we denoted $h_{i} = h_{i} (C,\nu)$ with $C\in  \left\{ 0.5, 1, 2, 10 \right\}$ and $\nu\in \left\{ \frac{1}{10} , \frac{1}{8} , \frac{1}{6} , \frac{1}{5} , \frac{1}{4} , \frac{1}{3} , \frac{1}{2} , 1 \right\}$, and we considered the cross validation criterion
$$ CV (C,\nu) = \frac{1}{n} \sum_{i=1}^{n} \left( Y_{i} - r_{n}^{[\ell],[-i]} (\mathcal{X}_{i}) \right)^{2}, $$
\noindent where $r_{n}^{[\ell],[-i]}$ represents the recursive estimator of $r$ using the $(n-1)$-sample that corresponds to the initial sample without the $i^{\text{th}}$ observation $(\mathcal{X}_{i},Y_{i})$ for $1\leq i \leq n$. Then, we selected the values $C_{CV}$, $\nu_{CV}$ of $C$ and $\nu$ that minimized $CV(C,\nu)$. Our estimator was $r_{n}^{[\ell]}$ using these selected values of $C$ and $\nu$.\\
Table 1 presents the mean and standard deviations of the prediction error over $500$ simulations for the optimal values of $C$ and $\nu$ with respect to the $CV$ criterion (these optimal values were $C_{CV}=1$ and $\nu_{CV} = 1/10$ for our estimator). Specifically, denoting $\widehat{Y}^{[j]} = r_{n}^{[\ell],[j]} (\chi^{[j]})$ the predicted value at the $j^{\text{th}}$ iteration of the simulation ($j = 1, \ldots, 500$) for a new curve $\chi^{[j]}$, we provided the mean ($MSPE$) and the standard deviations of the quantities $\left( \widehat{Y}^{[j]} - Y^{[j]} \right)^{2}$. The errors were computed for our estimator (label (1) in the table) and the estimator from \cite{Ferraty2006} (label (2) in the table), both of which were adapted using the \cite{RachdiVieu} procedure. From these results, we observed that the estimator from \cite{Ferraty2006} is slightly better than our estimator for the $MSPE$ criterion. As we will observe in (see Subsection \ref{computime}), the advantage of our estimator is its computational time. This section also examines the behavior of the prediction errors when the sample size increases. We used $n=100$, $n=200$ and $n=500$; as expected, the errors decreased when the sample size increased.
\begin{table}[h!]
\begin{center}
\begin{tabular}{cccc}
& $n=100$ & $n=200$ & $n=500$ \\
\hline
(1) & $0.3022$ & $0.2596$ & $0.1993$ \\
& $(0.6887)$ & $(0.6275)$ & $(0.5430)$ \\
\hline
(2) & $0.2794$ & $0.2143$ & $0.1368$ \\
& $(0.5512)$ & $(0.5055)$ & $(0.4208)$
\end{tabular}
\caption{Mean and standard deviation of the square prediction error computed via $500$ simulations for different values of $n$ with the optimal values of bandwidth provided by $C_{CV}$ and $\nu_{CV}$}
\end{center}
\end{table}
\subsection{Choice of the semi-norm}
In this simulation, the parameter $\ell$ was fixed equal to $0$ and we chose bandwidth $h_{i} = \max\limits_{i = 1, \ldots, n} \left\| \mathcal{X}_{i} - \chi \right\| i^{-1/10}$. Next, we aimed to compare the influence of the choice of the semi-norm by considering the following:

$\bullet$ The semi-norm $[PCA]$ based on the principal components analysis of the curves retained $q=3$ principal components; specifically,
$$ \left\| \mathcal{X}_{i} - \chi \right\|_{PCA} = \sqrt{\sum_{j=1}^{q} \langle \mathcal{X}_{i} - \chi , \nu_{j} \rangle^{2}}, $$
where $\langle . , . \rangle$ is the usual inner product of the space of square-integrable functions and $(\nu_{j})$ is the sequence of the eigenfunctions of the empirical covariance operator $\Gamma_{n}$ defined by $\Gamma_{n} u := \frac{1}{n} \sum_{i=1}^{n} \langle \X_{i} , u \rangle u$.

$\bullet$ The semi-norm $[FOU]$ was based on a decomposition of the curves in a Fourier basis, with $b=8$ basis functions; specifically,

$$ \left\| \mathcal{X}_{i} - \chi \right\|_{FOU} = \sqrt{\sum_{j=1}^{b} \left( a_{\mathcal{X}_{i},j} - a_{\chi,j} \right)^{2}}, $$

\noindent where $(a_{\mathcal{X}_{i},j})$ and $(a_{\chi,j})$ are the coefficient sequences of the respective Fourier approximations of curves $\mathcal{X}_{i}$ and $\chi$.

$\bullet$ The semi-norm $[DERIV]$ was based on a cubic spline approximation comparison of the second derivatives of the curves (with numerous interior knots $k=8$ for the cubic splines); specifically,

$$ \left\| \mathcal{X}_{i} - \chi \right\|_{DERIV} = \sqrt{\langle \widetilde{\mathcal{X}}_{i} - \widetilde{\chi} , \widetilde{\mathcal{X}}_{i} - \widetilde{\chi} \rangle}, $$

\noindent where $\widetilde{\mathcal{X}}_{i}$ and $\widetilde{\chi}$ are the spline approximations of the second derivatives of curves $\mathcal{X}_{i}$ and $\chi$.

$\bullet$ The semi-norm $[PLS]$ in which the data are projected on the space determined by a PLS regression on the curves applied $K=5$ PLS basis functions; specifically,

$$ \left\| \mathcal{X}_{i} - \chi \right\|_{PLS} = \sqrt{\sum_{j=1}^{K} \langle \mathcal{X}_{i} - \chi , p_{j} \rangle^{2}}, $$

\noindent where $(p_{j})$ is the sequence of PLS basis functions.

The results are provided in Table 2. For these simulated data, the semi-norms $[PCA]$ and $[PLS]$ show better results. As \cite{Ferraty2006} indicated, however, a universal norm does not exist that would overcome the others. The choice of the semi-norm depends on the treated data.

\begin{table}[h!]
\begin{center}
\begin{tabular}{ccccc}
norm & $[PCA]$ & $[FOU]$ & $[DERIV]$ & $[PLS]$ \\ \hline
MSPE & $0.3936$ & $0.4506$ & $0.4527$ & $0.3887$ \\
& $(1.5190)$ & $(1.5624)$ & $(1.5616)$ & $(1.5098)$
\end{tabular}
\caption{Mean and standard deviation of the square prediction error computed on $500$ simulations for different norms.}
\end{center}
\end{table}

\subsection{Choice of the parameter $\ell$}\label{choixl}
In this simulation, we used $h_{i} = \max\limits_{i = 1, \ldots, n} \left\| \X_{i} - \chi \right\| i^{-1/10}$ and the semi-norm based on the principal components analysis of the curves that retained $3$ principal components. The parameter $\ell$ varied into $\left\{ 0, \frac{1}{4}, \frac{1}{2},  \frac{3}{4}, \text{ and }1 \right\}$. The results are provided in Table 3. We showed that the values of the $MSPE$ criterion are close; thus, this parameter does not appear to significantly affect the quality of the prediction, even if the mean square error decreases with respect to $\ell$ as in the multivariate setting.

\begin{table}[h!]
\begin{center}
\begin{tabular}{cccccc}
 $\ell$ & $0$ & $0.25$ & $0.5$ & $0.75$ & $1$ \\ \hline
MSPE& $0.4054848$ & $0.4054814$ & $0.4054786$ & $0.4054764$ & $0.4054746$ \\
& $(1.372965)$ & $(1.372930)$ & $(1.372896)$ & $(1.372863)$ & $(1.372831)$
\end{tabular}
\caption{Mean and standard deviation of the square prediction error computed on $500$ simulations for different values of $\ell$.}
\end{center}
\end{table}

\subsection{Computational time} \label{computime}

This subsection highlights an important advantage of the recursive estimator compared with the initial one from \cite{Ferraty2006}), regarding the increase in computational time needed to predict a response when new values of the explanatory variable are sequentially added to the database. In fact, when a new observation $(\mathcal{X}_{n+1} , Y_{n+1})$ appears, the computation of the recursive estimator $r_{n+1}^{[\ell]}$ requests another iteration of the algorithm via its value computed with the sequence $\left( \mathcal{X}_{i} , Y_{i} \right)_{i = 1, \ldots, n}$, whereas the initial estimator must be recalculated using the whole sample $\left( \mathcal{X}_{i} , Y_{i} \right)_{i = 1, \ldots, n+1}$. The following illustrates the computation time difference between both estimators in these types of situations. From an initial sample $\left( \mathcal{X}_{i}, Y_{i} \right)_{i = 1, \ldots, n}$ with size $n=100$, we considered $N$ additional observations for different values of $N$. We compared the cumulative computational times to obtain the recursive and non-recursive estimators when adding these new observations. The characteristics of the computer we used to perform these computations were CPU: Duo $E4700$ $\, 2.60$ GHz, HD: $149$ Go, Memory: $3.23$ Go. The simulation was performed under the following conditions: the curves $\mathcal{X}_{1}, \ldots, \text{ and }\mathcal{X}_{n}$ as well as the new observations $\mathcal{X}_{n+1}, \ldots, \text{ and } \mathcal{X}_{n+N}$, which are standard Brownian motions on $[0,1]$ with $n = 100$ and $N \in \left\{ 1, 50, 100, 200, 500 \right\}$.  The semi-norm, bandwidth sequence and the parameter $\ell$ were chosen as in the previous cases.\\
The computational times are displayed in Table 4. Here, our estimator shows a clear advantage with regard to computational time compared with the estimator from \cite{Ferraty2006}.
\begin{table}[!h]
\begin{tabular}{cccccc}
$N$ & $1$ & $50$ & $100$ & $200$ & $500$ \\ \hline
comp. time for $r_{n+1}^{[\ell]}, \ldots, r_{n+N}^{[\ell]}$  & $0.125$ & $0.484$ & $0.859$ & $1.563$ & $3.656$ \\ \hline
comp. time for $r_{n+1}, \ldots, r_{n+N}$ & $0.047$ & $1.922$ & $5.594$ & $21.938$ & $152.719$
\end{tabular}
\caption{Cumulative computational times in seconds for the recursive and \cite{Ferraty2006} estimators when adding $N$ new observations for different values of $N$.}
\end{table}

\subsection{First real dataset example}

\begin{figure}[h!]
\includegraphics[width=12cm,height=8cm]{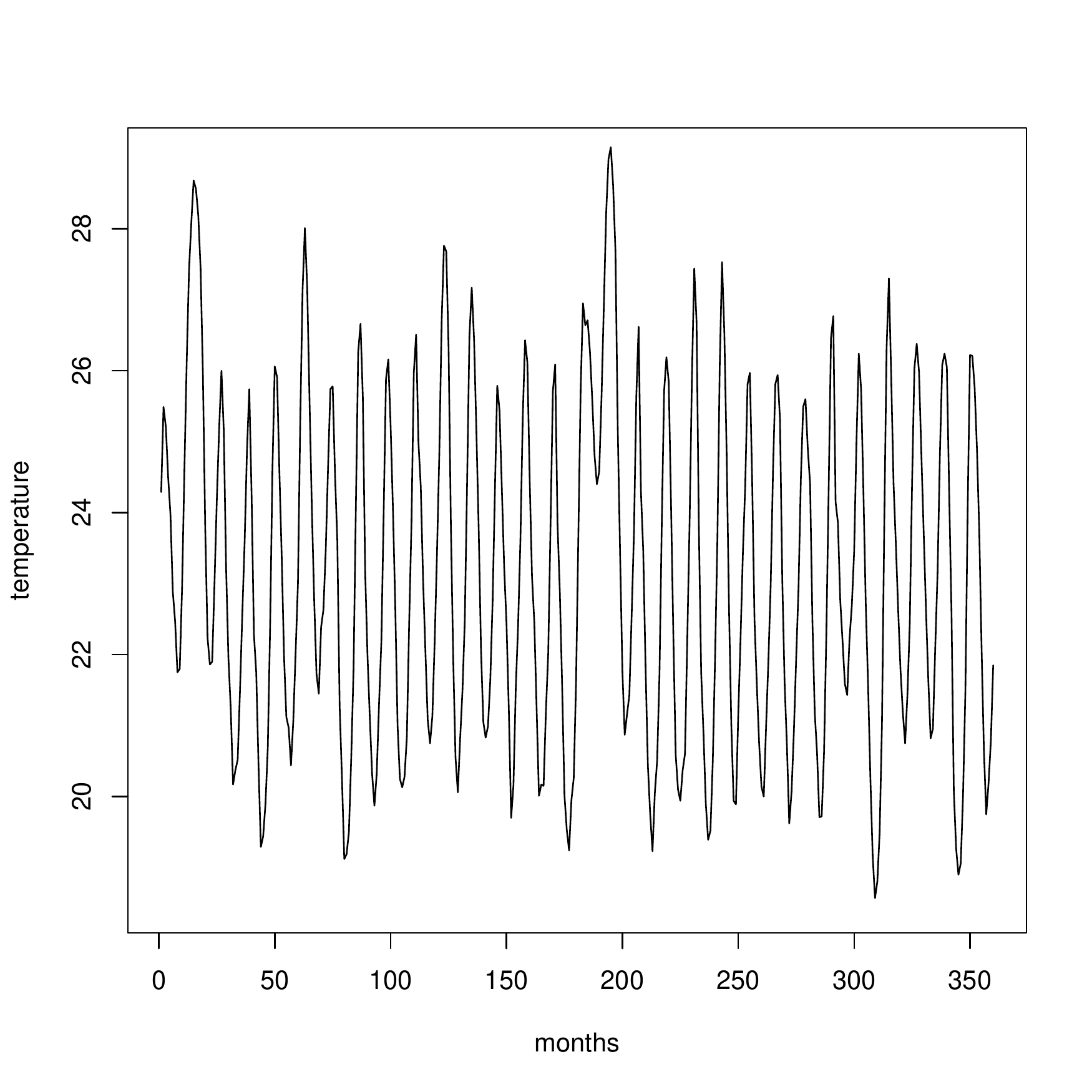}
\caption{El Ni$\tilde{\text{n}}$o monthly temperature time series from January 1982 to December 2011}
\label{fig1}
\end{figure}

\begin{figure}[h!]
\includegraphics[width=12cm,height=8cm]{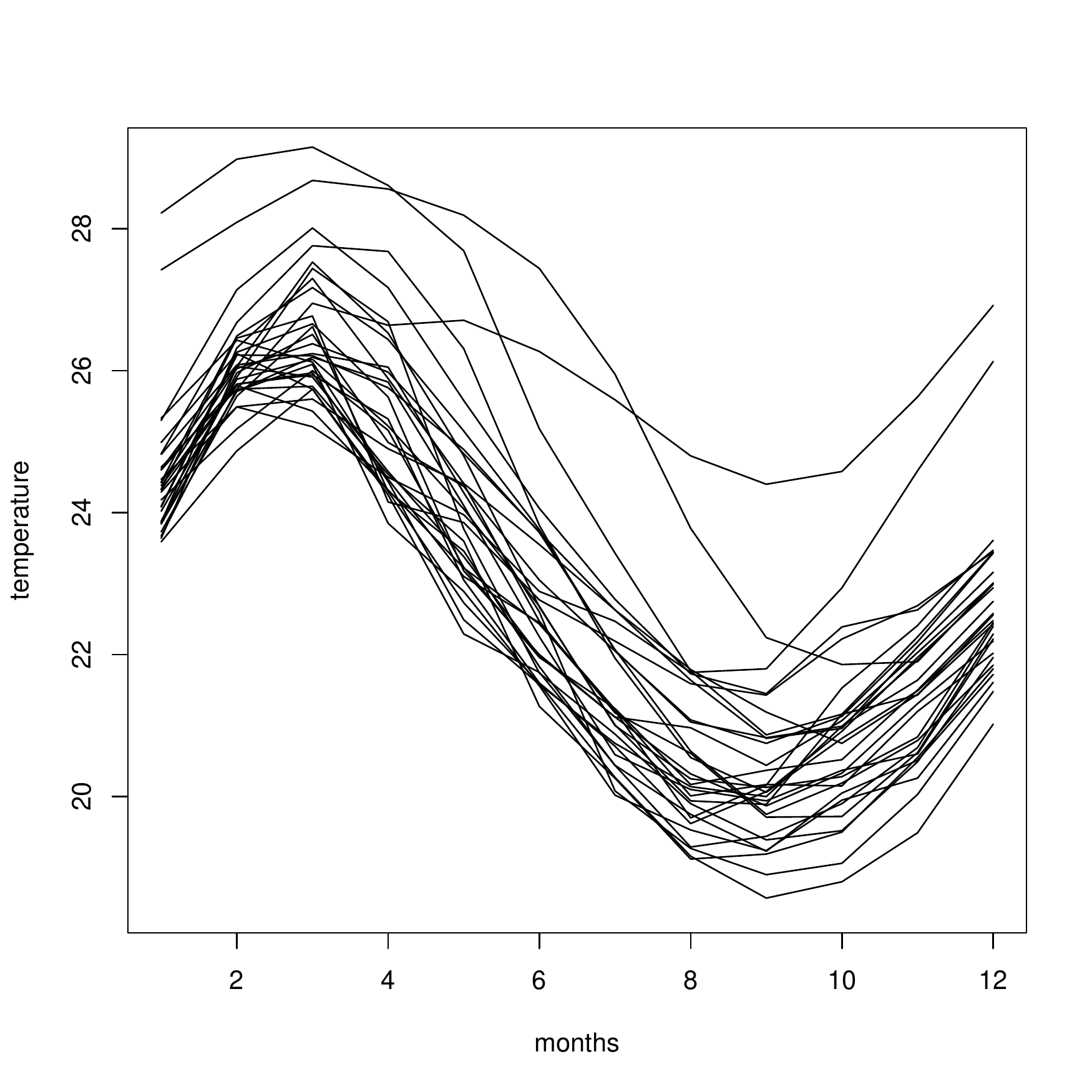}
\caption{El Ni$\tilde{\text{n}}$o annual temperature curves from 1982 to 2011}
\label{fig2}
\end{figure}

This subsection applies our estimator to a real dataset. Functional data are particularly adapted to studying a time series. We illustrate this fact using the El Ni$\tilde{\text{n}}$o time series\footnote{available online at http://www.math.univ-toulouse.fr/staph/npfda/} that provides the monthly sea surface temperature from January 1982 to December 2011 ($360$ months). These data are plotted in Figure \ref{fig1}. From this time series, we extracted the $30$ annual curves $\X_{1}, \ldots, \text{ and } \X_{30}$ from 1982 to 2011, discretized into $p = 12$ points. These curves are plotted in Figure \ref{fig2}. The variable of interest at month $j$ of year $i$ is the sea temperature $\X_{i+1}$ for month $j$; in other words, for $j = 1, \ldots, 12$ and for $i = 1, \ldots, 29,$ $ Y_{i}^{[j]} = \X_{i+1} (j). $

We predicted the values of $Y_{29}^{[1]}  \ldots, \text{ and } Y_{29}^{[12]}$ (in other words, the values of curve $\X_{30}$). The recursive estimator and the estimator from \cite{Ferraty2006} were computed by choosing the semi-norm, the bandwidth sequence and the parameter $\ell$ as in previous cases.

We analyzed the results by computing the mean square prediction error through 2011, which was provided by

$$ MSPE = \frac{1}{12} \sum_{j=1}^{12} \left( \widehat{Y}_{29}^{[j]} - Y_{29}^{[j]} \right)^{2}, $$

\noindent where $\widehat{Y}_{29}^{[j]}$ was computed with either the recursive estimator (result = $0.5719$) or the estimator from \cite{Ferraty2006} (result: $0.2823$). The estimator from \cite{Ferraty2006} again revealed its advantage with regard to prediction, whereas our estimator behaved well and had an advantage with regard to computational time (as highlighted in the previous subsection). The computational time (in seconds) needed for our estimator to predict twelve values (the final year) was $0.128$s, whereas the computational time for the \cite{Ferraty2006} estimator was $0.487$s.\\

To construct the confidence intervals described at the end of Section 2, we used the example from 2011.
The constants involved in Corollary \ref{normality2} were estimated by
\begin{eqnarray*}
\widehat M_1=\frac{1}{n}\sum_{i=1}^n\frac{1}{\widehat F(h_i)}K\left(\frac{\|\chi-\X_i\|}{h_i}\right)\text{ and }
\widehat M_2=\frac{1}{n}\sum_{i=1}^n\frac{1}{\widehat F(h_i)}K^2\left(\frac{\|\chi-\X_i\|}{h_i}\right).
\end{eqnarray*}
The consistency of $\widehat M_i, i=1,2$ is the same as the lines of the Lemmas \ref{Biais2}, and \ref{lem1}   proofs with the help of the Glivenko-Cantelli theorem. 
In addition, the estimation of the conditional variance $\sigma_\varepsilon^2$ was performed using a nonparametric kernel regression procedure:
\begin{eqnarray*}
\widehat\sigma_\varepsilon^2(\chi)=\frac{\sum\limits_{i=1}^nY_i^2K\left( \frac{\|{\chi}-{\cal X}_i\|}{h_i}\right)}
{\sum\limits_{i=1}^nK\left( \frac{\|{\chi}-{\cal X}_i\|}{h_i}\right)}-\left(r_n^{[0]}(\chi)\right)^2.
\end{eqnarray*}
The different values of $\widehat\sigma_{\varepsilon}^{2}$ for each month of last year are shown in Table 5.

\begin{table}[h!]
\begin{center}
\begin{tabular}{ccccccc}
month & $1$ & $2$ & $3$ & $4$ & $5$ & $6$ \\
\hline
$\widehat{\sigma}_{\varepsilon}^{2}$ & $1.00551$ & $0.60462$ & $0.81971$ & $1.60694$ & $2.07112$ & $2.02370$ \\
\hline
month & $7$ & $8$ & $9$ & $10$ & $11$ & $12$ \\
\hline
$\widehat{\sigma}_{\varepsilon}^{2}$ & $2.01471$ & $1.65082$ & $1.25991$ & $1.24542$ & $1.27177$ & $1.08482$
\end{tabular}
\caption{Recursive estimations of $\sigma_{\varepsilon}^{2}$ for each month on the last year.}
\end{center}
\end{table}

\begin{figure}[h!]
\includegraphics[width=12cm,height=8cm]{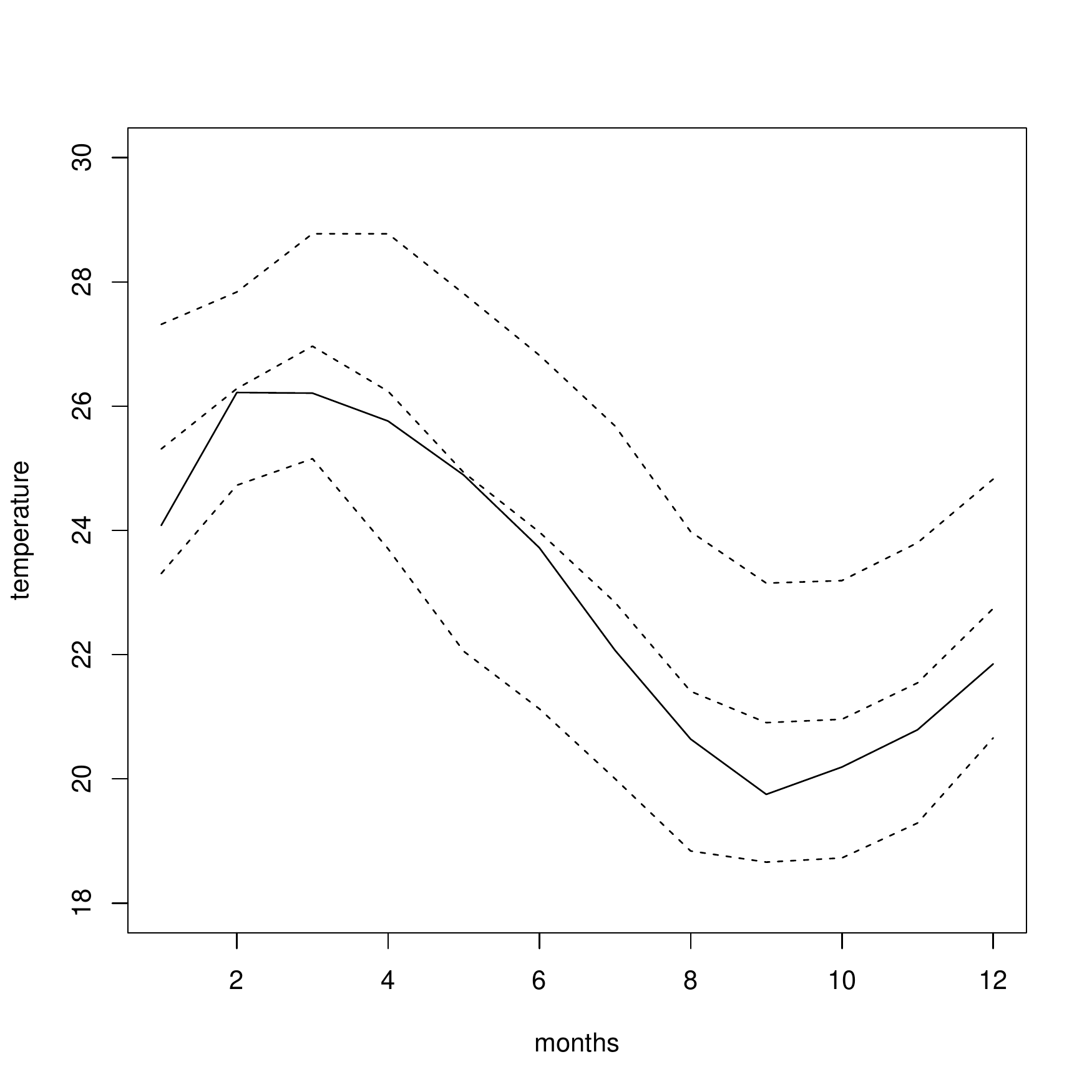}
\caption{El Ni$\tilde{\text{n}}$o true and predicted temperature curves for 2011. The solid line denotes the true curve. The dashed lines represent the predicted curve and the $95 \%$ confidence intervals with the recursive estimator.}
\label{fig3}
\end{figure}

\noindent The estimations of the other quantities in the confidence intervals were $\widehat{\beta}_{[1]} = 1.04966$, $\widehat{M}_{1} = 0.73367$ and $\widehat{M}_{2} = 0.62716$.\\
We compared the confidence intervals for 2011 in Table 6 with the recursive estimator (label (1) in the table) and the estimator from \cite{Ferraty2006} (label (2) in the table). We noticed that the confidence interval with the estimator from \cite{Ferraty2006} was still better, mainly due to the estimation of $\sigma_{\varepsilon}$ and the terms $M_{1}$ and $M_{2}$ (which were estimated by the respective non-recursive counterparts of $\widehat{\sigma}_{\varepsilon}^{2}$ , and $\widehat{M}_{1}$, and $\widehat{M}_{2}$ was obtained by replacing $h_i$ with $h_n$ in the non-recursive case), whereas the term $\widehat{\beta}_{[1]}$ decreased the width of the confidence interval in the recursive case.\\
The corresponding true curve, predicted curve and $95 \%$ confidence intervals over the year 2011 are plotted in Figure \ref{fig3}.
\begin{center}
\begin{table}[h!]
 \begin{tabular}{ccccccc}
month & $1$ & $2$ & $3$ & $4$  \\
\hline
(1) & $[23.24;27.38]$ & $[24.68;27.89]$ & $[25.10;28.83]$ & $[23.62;28.86]$ &  \\
\hline
(2) & $[22.80;26.24]$ & $[24.83;27.49]$ & $[25.07;28.17]$ & $[23.42;27.76]$ & \\
\hline
month & $5$ & $6$ &$7$ & $8$ &\\
\hline
(1) & $[21.96;27.90]$ & $[21.04;26.91]$& $[19.91;25.77]$ & $[18.76;24.06]$ &\\
\hline
(2) & $[21.80;26.74]$ & $[20.44;25.32]$& $[19.18;24.05]$ & $[18.42;22.83]$ \\
\hline
month &  $9$ & $10$ & $11$ & $12$ \\
\hline
(1) & $[18.59;23.22]$ & $[18.65;23.26]$ & $[19.22;23.87]$ & $[20.59;24.89]$ \\
\hline
(2) & $[18.37;22.22]$ & $[18.81;22.63]$ & $[19.54;23.40]$ & $[20.88;24.45]$\\
\hline
\end{tabular}
\caption{Confidence intervals for 2011 with the recursive estimator (label (1)) and the estimator from \cite{Ferraty2006} (label (2)).}
\end{table}
\end{center}
\subsection{Second real dataset example}
This subsection highlights the performances of our estimator with regard to computational time gain.  The dataset we used consisted of information collected by ORAMIP\footnote{Observatoire R\'egional de l'Air en Midi-Pyr\'en\'ees}, a French organization that studies air quality. We disposed of a sample of $474$ daily ozone pollution measurements curves (in $\mu$g/m$^{3}$). The variable $Y_{i}$ of interest was the daily maximum of ozone. The ozone curve the day before (from $6:00$ pm to $5:00$ pm the day after) was used as a functional explicative variable $\X_{i}$. Specifically, each ${\X_{i}}$ was observed at p = 24 equidistant points that corresponded to hourly measurements. The sample was divided into a learning sample of $332$ and a test sample of $142$. The ozone curve from the learning sample is plotted in Figure 4.
\begin{figure}[h!]
\includegraphics[width=12cm,height=8cm]{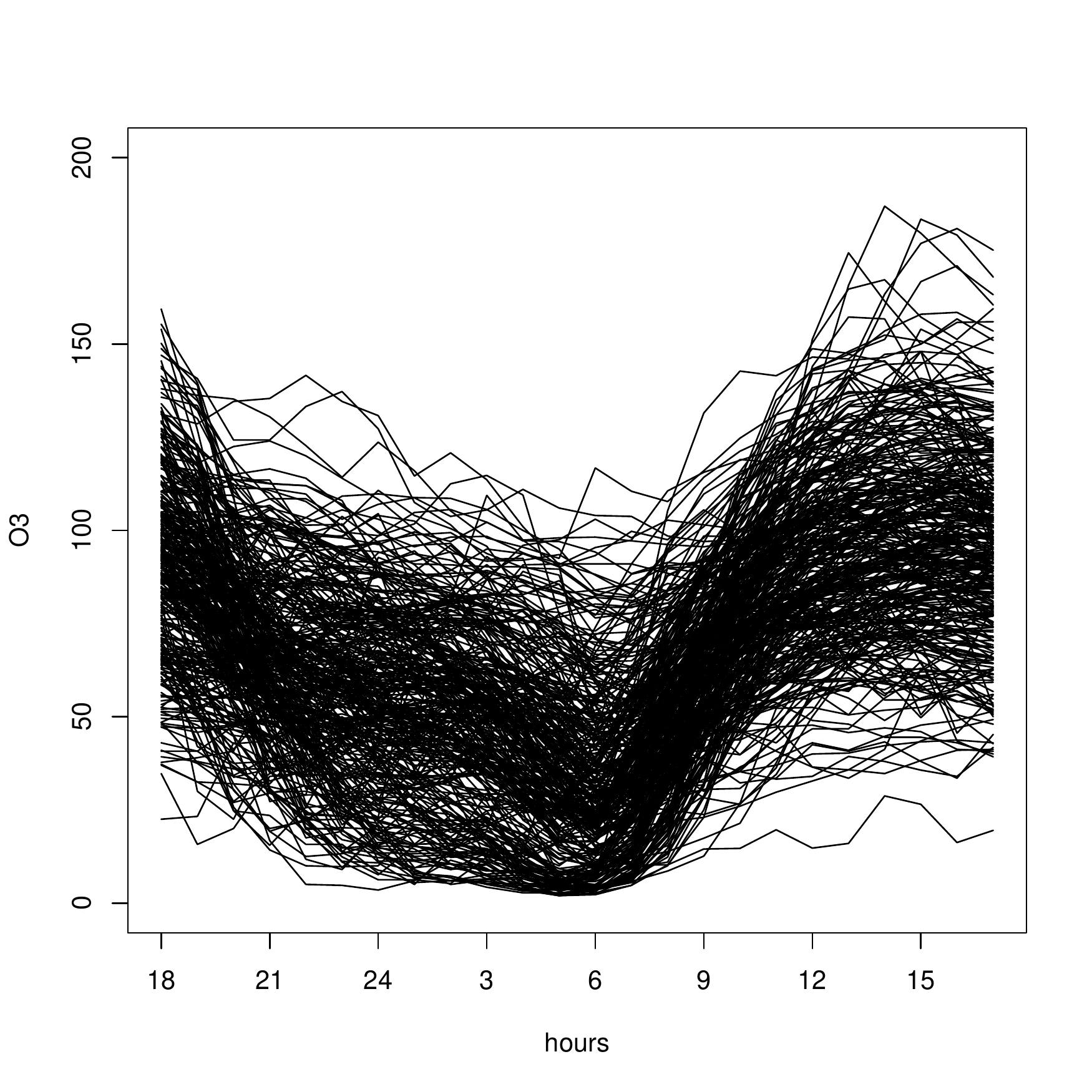}
\caption{Daily ozone curves for the learning sample}
\end{figure}
Every value of the test sample was predicted using our method and the method from \cite{Ferraty2006}. The learning sample and the curves of the test sample were predicted as if they had arrived in real time. Again, the estimator from \cite{Ferraty2006} showed its advantage over our estimator with regard to mean prediction error ($25.797$ vs. $26.171$). However, our estimator performed best when predicting several values and using recursivity, especially when many values need to be predicted (which was the case here, given that we predicted $142$ values). The computational time for our estimator (in seconds) was $2.025$, whereas it was $25.478$ for the Ferraty and Vieu method; however, we are conscious that this advantage was relative in this case because we dealt with daily measurements. In any case, the choice is up to the user: a better prediction error or a faster computational time.

\section{Proofs}\label{proof}
Throughout  the proofs, we denote $\gamma_i$ a sequence of real numbers going to zero as $i$ goes to $\infty$.
The  kernel estimate $r_n^{[\ell]}$ writes
\begin{eqnarray*}
r_n^{[\ell]}(\chi)&=&\frac{\varphi_n^{[\ell]}({\chi})}{f_n^{[\ell]}({\chi})},
\end{eqnarray*}
where  $\varphi_n^{[\ell]}$ and $f_n^{[\ell]}$ are  defined in (\ref{phin}).
\subsection{Proof of Theorem \ref{Bias_Var}}
To prove the first assertion of Theorem \ref{Bias_Var}, let  us use the following decomposition
\begin{eqnarray*}
\text{E}\left[r_n^{[\ell]}({\chi})\right]&=&\frac{\text{E}\left[\varphi_n^{[\ell]}({\chi})\right]}{\text{E}\left[f_n^{[\ell]}({\chi})\right]}-\frac{\text{E}\left\{\varphi_n^{[\ell]}({\chi})\left[f_n^{[\ell]}({\chi})-\text{E}f_n^{[\ell]}({\chi})\right]\right\}}{\left\{\text{E}\left[f_n^{[\ell]}({\chi})\right]\right\}^2}\\ &&+\frac{\text{E}\left\{r_n^{[\ell]}({\chi})\left[f_n^{[\ell]}({\chi})-\text{E}f_n^{[\ell]}({\chi})\right]^2\right\}}{\left\{\text{E}\left[f_n^{[\ell]}({\chi})\right]\right\}^2}.
\end{eqnarray*}
The first part of Theorem \ref{Bias_Var} is then a direct consequence of Lemmas \ref{Biais1}-\ref{Var} below.
\begin{lem}\label{Biais1}
Under assumptions \textbf{H1}-\textbf{H4}, we have
\begin{eqnarray*}
\frac{\text{E}\left[\varphi_n^{[\ell]}({\chi})\right]}{\text{E}\left[f_n^{[\ell]}({\chi})\right]}-r(\chi)=h_n\varphi'(0)\frac{\alpha_{[\ell]}}{\beta_{[1-\ell]}}\frac{M_0}{M_1}\left[1+o(1)\right].
\end{eqnarray*}
\end{lem}
\begin{lem}\label{termesA}
Under assumptions \textbf{H1}-\textbf{H4}, we have
\begin{eqnarray*}
\text{E}\left\{\varphi_n^{[\ell]}({\chi})\left[f_n^{[\ell]}({\chi})-\text{E}f_n^{[\ell]}({\chi})\right]\right\}&=& O\left[\frac{1}{nF(h_n)}\right],\\
\text{E}\left\{r_n^{[\ell]}({\chi})\left[f_n^{[\ell]}({\chi})-\text{E}f_n^{[\ell]}({\chi})\right]^2\right\}&=&O\left[\frac{1}{nF(h_n)}\right].
\end{eqnarray*}
\end{lem}

\begin{lem}\label{Biais2}
Under assumptions \textbf{H1}-\textbf{H4}, we have
\begin{eqnarray*}
\text{E}\left(f_n^{[\ell]}({\chi})\right)=M_1\left[1+o(1)\right] & \text{ and }&
\text{E}\left(\varphi_n^{[\ell]}({\chi})\right)=r(\chi)M_1\left[1+o(1)\right].
\end{eqnarray*}
\end{lem}
\noindent The  study of  the variance term in Theorem \ref{Bias_Var} is based  on  the following decomposition  which can be found in  \cite{Collomb1976}.
\begin{eqnarray}\label{egalvar}
\text{Var}\left[r_n^{[\ell]}(\chi)\right]&=& \frac{\text{Var}\left[\varphi_n^{[\ell]}({\chi})\right]}{\left\{\text{E}\left[f_n^{[\ell]}({\chi})\right]\right\}^2}-4\frac{\Eb\left[\varphi_n^{[\ell]}({\chi})\right]\text{Cov}\left[f_n^{[\ell]}({\chi}),\varphi_n^{[\ell]}({\chi})\right]}{\left\{\Eb\left[f_n^{[\ell]}({\chi})\right]\right\}^3} \nonumber\\\mbox{}&& +3\text{Var}\left[f_n^{[\ell]}({\chi})\right]\frac{\left\{\Eb\left[\varphi_n^{[\ell]}({\chi})\right]\right\}^2}{\left\{\Eb\left[f_n^{[\ell]}({\chi})\right]\right\}^4}+o\left[\frac{1}{nF(h_n)}\right].
\end{eqnarray}
The second assertion of Theorem \ref{Bias_Var} follows from  \eqref{egalvar} and Lemma \ref{Var} below. \hfill $\square$
\begin{lem}\label{Var}
Under assumptions \textbf{H1}-\textbf{H4}, we have
\begin{eqnarray*}
\text{Var}\left[f_n^{[\ell]}({\chi})\right]&=&\frac{\beta_{[1-2\ell]}}{\beta_{[1-\ell]}^2}M_2\frac{1}{nF(h_n)}\left[1+o(1)\right].\\
\text{Var}\left[\varphi_n^{[\ell]}({\chi})\right]&=&\frac{\beta_{[1-2\ell]}}{\beta_{[1-\ell]}^2}\left[r^2(\chi)+\sigma_\epsilon^2(\chi)\right]M_2\frac{1}{nF(h_n)}\left[1+o(1)\right].\\
\text{Cov}\left[f_n^{[\ell]}({\chi}),\varphi_n^{[\ell]}({\chi})\right]&=&\frac{\beta_{[1-2\ell]}}{\beta_{[1-\ell]}^2}r(\chi)M_2\frac{1}{nF(h_n)}\left[1+o(1)\right].
\end{eqnarray*}
\end{lem}
Now let us  prove Lemmas \ref{Biais1}-\ref{Var}.
\subsubsection{Proof of Lemma \ref{Biais1}}

Observe that
\begin{eqnarray*}
\frac{\text{E}\left[\varphi_n^{[\ell]}({\chi})\right]}{\text{E}\left[f_n^{[\ell]}({\chi})\right]}-r(\chi)&=&\dfrac{\sum\limits_{i=1}^n\frac{1}{F(h_i)^\ell}\text{E}\left[\left(Y_i-r(\chi)\right)K\left(\frac{\|\chi-{\cal X}_i\|}{h_i}\right)\right]}{\sum\limits_{i=1}^n\frac{1}{F(h_i)^\ell}\text{E}\left[K\left(\frac{\|\chi-{\cal X}_i\|}{h_i}\right)\right]}.
\end{eqnarray*}
Obviously 
\begin{eqnarray*}
\Eb\left[\left(Y_i-r(\chi)\right)K\left(\frac{\|\chi-{\cal X}_i\|}{h_i}\right)\right]&=&\Eb\left[\left(r(\X)-r(\chi)\right)K\left(\frac{\|\X-\chi\|}{h_i}\right)\right]\\&=&
\Eb\left[\varphi\left(\|\X-\chi\|\right)K\left(\frac{\|\X-\chi\|}{h_i}\right)\right]\\&=& \int_0^1\varphi(h_it)K(t)d\Pb^{\|\X-\chi\|/h_i}(t).
\end{eqnarray*}
Therefore, a Taylor's expansion of $\varphi$ around 0 ensures that
\begin{eqnarray*}
\Eb\left[\left(Y_i-r(\chi)\right)K\left(\frac{\|\chi-{\cal X}_i\|}{h_i}\right)\right]&=&h_i\varphi'(0)\int_0^1 tK(t)d\Pb^{\|\X-\chi\|/h_i}(t)+o(h_i).
\end{eqnarray*}
From the proof of Lemma 2 in \cite{Ferraty2007},  {\bf H2} and Fubini's Theorem imply that
\begin{eqnarray}\label{equation1}
\int_0^1 tK(t)d\Pb^{\|\X-\chi\|/h_i}(t)&=&F(h_i)\left[K(1)-\int_0^1(sK(s))'\tau_{h_i}(s)ds\right],
\end{eqnarray}
and
\begin{eqnarray}\label{equation2}
\Eb K\left(\frac{\|\X-\chi\|}{h_i}\right)=\int_0^{h_i}K\left(\frac{t}{h_i}\right)d\Pb^{\|\X-\chi\|}(t)=F(h_i)\left[K(1)-\int_0^1K'(s)\tau_{h_i}(s)ds\right].
\end{eqnarray}
Also,  combining \eqref{equation1} and \eqref{equation2}, we have
\begin{eqnarray*}
\frac{\text{E}\left[\varphi_n^{[\ell]}({\chi})\right]}{\text{E}\left[f_n^{[\ell]}({\chi})\right]}-r(\chi)&=&\frac{\sum\limits_{i=1}^nh_iF(h_i)^{1-\ell}\left\{\varphi'(0)\left[K(1)-\int_0^1(sK(s))'\tau_{h_i}(s)ds\right]+\gamma_i\right\}}{\sum\limits_{i=1}^nF(h_i)^{1-\ell}\left[K(1)-\int_0^1K'(s)\tau_{h_i}(s)ds\right]}\\&:=&\frac{D_1}{D_2}.
\end{eqnarray*}
Finally, by virtue of {\bf H3} we get from  Toeplitz's Lemma (see \cite{Masry86}) that
\begin{eqnarray*}
\frac{D_1}{nh_n{F(h_n)}^{1-\ell}}=\alpha_{[\ell]}\varphi'(0) M_0[1+o(1)], \ \
\frac{D_2}{n{F(h_n)}^{1-\ell}}=\beta_{[1-\ell]}M_1[1+o(1)],
\end{eqnarray*}
and Lemma \ref{Biais1} follows. \hfill $\square$
\subsubsection{Proof of Lemma \ref{Biais2}}
Equality \eqref{equation2}  allows to  write 
\begin{eqnarray*}
\Eb\left[f_n^{[\ell]}({\chi})\right]&=&\frac{1}{\sum\limits_{i=1}^n{F(h_i)}^{1-\ell}}\sum\limits_{i=1}^n\frac{1}{F(h_i)^\ell}\Eb\left[K\left(\frac{\|\chi-\X_i\|}{h_i}\right)\right]\\&=& 
\frac{\sum\limits_{i=1}^n \frac{F(h_i)^{1-\ell}}{nF(h_n)^{1-\ell}}\left[K(1)-\int_0^1K'(s)\tau_{h_i}(s)ds\right]}{\displaystyle B_{n,1-\ell}}=M_1[1+o(1)],
\end{eqnarray*}
where the last equality follows from assumptions \textbf{H3}, \textbf{H4} and Toeplitz's Lemma. Now, conditioning on $\X$, we have
\begin{eqnarray*}
\Eb\left[Y_iK\left(\frac{\|\chi-\X_i\|}{h_i}\right)\right]&=&\Eb\left\{\left[r(\X)-r(\chi)+r(\chi)\right]K\left(\frac{\|\chi-\X_i\|}{h_i}\right)\right\}=:A_i+B_i,
\end{eqnarray*}
where
$$ A_i:=\Eb\left\{\left[r(\X)-r(\chi)\right]K\left(\frac{\|\chi-\X_i\|}{h_i}\right)\right\} \leq \sup_{\chi'\in \mathcal{B}(\chi,h_i)}\left|r(\chi')-r(\chi)\right| \Eb K\left(\frac{\|\chi-\X_i\|}{h_i}\right), $$
and  $ B_i:=r(\chi)\Eb K\left(\frac{\|\chi-\X_i\|}{h_i}\right). $
Since $r$ is continuous (\textbf{H1}), then
\begin{eqnarray}\label{esperance}
\Eb\left[Y_iK\left(\frac{\|\chi-\X_i\|}{h_i}\right)\right]= \left[r(\chi)+\gamma_i\right]\Eb K\left(\frac{\|\chi-\X_i\|}{h_i}\right)=F(h_i)M_1\left[r(\chi)+\gamma_i\right].
\end{eqnarray}
We deduce from \eqref{esperance}, with the help of assumptions \textbf{H3} and \textbf{H4}, and applying again Toeplitz's Lemma, that
\begin{eqnarray*}
\Eb\left[\varphi_n^{[\ell]}({\chi})\right]&=& \frac{1}{\sum\limits_{i=1}^n{F(h_i)}^{1-\ell}}\sum_{i=1}^n\frac{1}{F(h_i)^\ell}\Eb\left[Y_iK\left(\frac{\|\chi-\X_i\|}{h_i}\right)\right] =r(\chi)M_1\left[1+o(1)\right],
\end{eqnarray*}
that proves Lemma \ref{Biais2}. \hfill$\square$
\subsubsection{Proof of Lemma \ref{Var}}
First, notice  that as in \eqref{equation2}, we have
\begin{eqnarray}\label{M2}
\Eb\left[K^2\left(\frac{\|\chi-\X\|}{h_i}\right)\right]&=&F(h_i)\left[K^2(1)-\int_0^1(K^2)'(s)\tau_{h_i}(s)ds\right].
\end{eqnarray}
Relation  \eqref{equation2} and assumption \textbf{H3} ensure that
\begin{eqnarray*}
\Eb^2\left[K\left(\frac{\|\chi-\X\|}{h_i}\right)\right]&=&O\left[F(h_i)^2\right],
\end{eqnarray*}
thus
\begin{eqnarray*}
\text{Var}\left[K\left(\frac{\|\chi-\X\|}{h_i}\right)\right]=M_2F(h_i)\left[1+\gamma_i\right].
\end{eqnarray*}
It follows  that
\begin{eqnarray*}
\text{Var}\left[f_n^{[\ell]}({\chi})\right]&=& \frac{1}{\left(\sum\limits_{i=1}^nF(h_i)^{1-\ell}\right)^2}\sum_{i=1}^nF(h_i)^{1-2\ell}M_2\left[1+\gamma_i\right]\\
&=&\frac{\beta_{[1-2\ell]}}{\beta_{[1-\ell]}^2}\frac{1}{nF(h_n)}M_2\left[1+o(1)\right],
\end{eqnarray*}
and the first step of Lemma \ref{Var} follows. In a similar manner, for the second step, let us  write
\begin{eqnarray*}
\text{Var}\left[\varphi_n^{[\ell]}({\chi})\right]&=&\frac{1}{\left(\sum\limits_{i=1}^nF(h_i)^{1-\ell}\right)^2}\sum_{i=1}^nF(h_i)^{-2\ell}\text{Var}\left[Y_iK\left(\frac{\|\chi-\X_i\|}{h_i}\right)\right].
\end{eqnarray*}
Next, one obtains by conditioning on $\X$,
\begin{eqnarray*}
\Eb\left[ Y_i^2K^2\left(\frac{\|\chi-\X_i\|}{h_i}\right)\right]&=& \Eb\left[r^2(\X)K^2\left(\frac{\|\chi-\X_i\|}{h_i}\right)\right]+\Eb\left[\sigma_\varepsilon^2(\X)K^2\left(\frac{\|\chi-\X_i\|}{h_i}\right)\right].
\end{eqnarray*}
Assumption \textbf{H1} and \eqref{M2} ensure that
\begin{eqnarray*}
\Eb\left[Y_i^2K^2\left(\frac{\|\chi-\X_i\|}{h_i}\right)\right]&=&\left[r^2(\chi)+\sigma_\varepsilon^2(\chi)\right]\Eb\left[K^2\left(\frac{\|\chi-\X_i\|}{h_i}\right)\right]\left[1+\gamma_i\right]\\&=& \left[r^2(\chi)+\sigma_\varepsilon^2(\chi)\right]M_2F(h_i)\left[1+\gamma_i\right],
\end{eqnarray*}
and then from Toeplitz's Lemma, with \textbf{H3} and \textbf{H4}, it follows that
\begin{eqnarray*}
\text{Var}\left[\varphi_n^{[\ell]}({\chi})\right]&=&\frac{1}{\left(\sum\limits_{i=1}^nF(h_i)^{1-\ell}\right)^2}\sum_{i=1}^nF(h_i)^{1-2\ell}\left[r^2(\chi)+\sigma_\varepsilon^2(\chi)\right]M_2\left[1+\gamma_i\right]
\\&=&\frac{\beta_{[1-2\ell]}}{\beta_{[1-\ell]}^2}\left[r^2(\chi)+\sigma_\varepsilon^2(\chi)\right]M_2\frac{1}{nF(h_n)}\left[1+o(1)\right],
\end{eqnarray*}
which proves the second assertion of Lemma \ref{Var}. The covariance term  writes
$$\begin{array}{rl}
\text{Cov}\left[f_n^{[\ell]}({\chi}),\varphi_n^{[\ell]}({\chi})\right]=&\frac{1}{\left(\sum\limits_{i=1}^nF(h_i)^{1-\ell}\right)^2}\Bigg\{\Eb\left[\sum\limits_{i=1}^n\sum\limits_{j=1}^n\frac{Y_iK\left(\frac{\|\chi-\X_i\|}{h_i}\right)K\left(\frac{\|\chi-\X_j\|}{h_j}\right)}{F(h_i)^{\ell}F(h_j)^{\ell}}\right]\\
 -&\sum\limits_{i=1}^n\frac{\Eb\left[Y_iK\left(\frac{\|\chi-\X_i\|}{h_i}\right)\right]}{F(h_i)^{\ell}}\sum\limits_{j=1}^n\frac{\Eb K\left(\frac{\|\chi-\X_j\|}{h_j}\right)}{F(h_j)^{\ell}}\Bigg\}\\
=&\frac{1}{\left(\sum\limits_{i=1}^nF(h_i)^{1-\ell}\right)^2}\sum\limits_{i=1}^n\frac{\Eb\left[Y_iK^2\left(\frac{\|\chi-\X_i\|}{h_i}\right)\right]}{F(h_i)^{2\ell}}\\
 -&\frac{1}{\left(\sum\limits_{i=1}^nF(h_i)^{1-\ell}\right)^2}\sum\limits_{i=1}^n\frac{\Eb\left[Y_iK\left(\frac{\|\chi-\X_i\|}{h_i}\right)\right]\Eb K\left(\frac{\|\chi-\X_i\|}{h_i}\right)}{F(h_i)^{2\ell}}:=I-II.
\end{array}$$
Notice that by  \eqref{equation1} and \eqref{esperance}, one may write 
\begin{eqnarray*}
II&=&O\left[ \frac{1}{n} \left( B_{n,1-\ell}\right)^{-2}B_{n,2(1-\ell)}\right]=O\left(\frac{1}{n F(h_n)}\right).
\end{eqnarray*}
Next, from assumption \textbf{H1} and conditioning on $\X$, we have
\begin{eqnarray*}
\Eb\left[Y_iK^2\left(\frac{\|\chi-\X_i\|}{h_i}\right)\right]=
M_2F(h_i)\left[r(\chi)+\gamma_i\right].
\end{eqnarray*}
It follows that
\begin{eqnarray*}
I= \frac{\left(\displaystyle B_{n,1-\ell}\right)^{-2}}{nF(h_n)}\sum_{i=1}^n\frac{F(h_i)^{1-2\ell}}{nF(h_n)^{1-2\ell}}M_2r(\chi)\left[1+\gamma_i\right],
\end{eqnarray*}
and the third assertion  of Lemma \ref{Var} follows again by applying Toeplitz's Lemma.\hfill $\square$
\subsubsection{Proof of Lemma \ref{termesA}}
Lemma \ref{termesA} is a direct consequence of Lemmas \ref{Biais2} and \ref{Var}.\hfill $\square$
\subsection{Proof of Theorem \ref{Cv_ps_rnl}}
Let us consider  the following decomposition
\begin{eqnarray}\label{egalite1}
 r_n^{[\ell]}(\chi)-r(\chi)=\frac{\tilde{\varphi}_n^{[\ell]}(\chi)-r(\chi)f_n^{[\ell]}(\chi)}{f_n^{[\ell]}(\chi)}+\frac{\varphi_n^{[\ell]}(\chi)-
\tilde{\varphi}_n^{[\ell]}(\chi)}{f_n^{[\ell]}(\chi)},
\end{eqnarray}
where $\tilde{\varphi}_n^{[\ell]}(\chi)$ is a truncated version of  $\varphi_n^{[\ell]}(\chi)$ defined by
 \begin{equation}\label{phi_nl_tilde}\tilde{\varphi}_n^{[\ell]}({\chi})=\frac{1}{\sum\limits_{i=1}^nF(h_i)^{1-\ell}}\sum\limits_{i=1}^n\dfrac{Y_i}{F(h_i)^{\ell}}\mathds{1}_{\left\{\left|Y_i\right|\leq b_n\right\}}K\left( \dfrac{\|{\chi}-{\cal X}_i\|}{h_i}\right),\end{equation}
$b_n$ being a sequence of real numbers which goes to $+\infty$ as $n\rightarrow\infty.$
Next, for any $\varepsilon>0$, we have for the residual term of (\ref{egalite1})
$$P\left\{\left|\varphi_n^{[\ell]}(\chi)-\tilde{\varphi}_n^{[\ell]}(\chi)\right|>\varepsilon\left[\frac{ \ln\ln n}{ nF(h_n)}\right]^\frac{1}{2}\right\}
 \leq P\left(\bigcup\limits_{i=1}^n\left\{|Y_i|>b_n\right\}\right)
  \leq\text{E}\left[e^{\lambda|Y|^\mu
 }\right]n^{1-\lambda\delta},$$
 where the last inequality follows by setting $b_n=(\delta\ln n)^\frac{1}{\mu},$  with the help of Markov's inequality.  Assumption ${\bf H5}$ ensures that for any $\varepsilon>0$, $$\sum_{n=1}^\infty P\left\{\left|\varphi_n^{[\ell]}(\chi)-\tilde{\varphi}_n^{[\ell]}(\chi)\right|>\varepsilon\left[\frac{ \ln\ln n}{ nF(h_n)}\right]^\frac{1}{2}\right\}<\infty \text{ if  $\delta>\frac{2}{\lambda}$} ,$$   which together with the Borel-Cantelli Lemma   imply that 
 \begin{eqnarray}\label{conv_residual}
 \left[\frac{ nF(h_n)}{\ln\ln n}\right]^{1/2}\left|\varphi_n^{[\ell]}(\chi)-\tilde{\varphi}_n^{[\ell]}(\chi)\right|\rightarrow 0 \text{ a.s,} \text{ as  }n\rightarrow\infty.
 \end{eqnarray}
The main term in \eqref{egalite1} writes 
 \begin{eqnarray}
 \label{decomposition11}
 \tilde{\varphi}_n^{[\ell]}(\chi)-r(\chi)f_n^{[\ell]}(\chi) &=&\left\{\tilde{\varphi}_n^{[\ell]}(\chi)-r(\chi)f_n^{[\ell]}(\chi)-\text{E}\left[\tilde{
\varphi}_n^{[\ell]}(\chi)-r(\chi)f_n^{[\ell]}(\chi)\right]\right\}\nonumber\\*
&&+\left\{ \text{E}\left[\tilde{\varphi}_n^{[\ell]}(\chi)-r(\chi)f_n^{[\ell]}(\chi)\right]\right\}:= N_1+N_2.
\end{eqnarray}
Theorem \ref{Cv_ps_rnl}  will be proved if  Lemmas \ref{lem1} and \ref{lem3} below are established. Indeed, from Lemma \ref{Biais2} we have $\text{E}\left(f_n^{[\ell]}({\chi})\right)=M_1\left[1+o(1)\right]$ and it can be shown (following the same lines of the proof of Lemma \ref{lem1}) that $$f_n^{[\ell]}(\chi)-\Eb f_n^{[\ell]}(\chi)=O\left(\sqrt{\frac{\ln\ln n}{nF(h_n)}}\right) \text{a.s.}$$  \hfill  $\square$
\begin{lem} \label{lem1} Under assumptions ${\bf H1-H6}$, we have
\begin{eqnarray*}
\label{I1reg}\varlimsup_{n\rightarrow\infty}\left[\frac{nF(h_n)}{\ln\ln~ n}\right]^{1/2}
N_1=\frac{\left[2\beta_{[1-2\ell]}\sigma^2_\varepsilon(\chi)M_2\right]^{1/2}}{\beta_{[1-\ell]}}\text{ a.s. }
\end{eqnarray*}
\end{lem}

{\lem \label{lem3}
Assume that ${\bf H1}-{\bf H5}$ hold. If $ \displaystyle\lim\limits_{n\rightarrow+\infty}nh_n^{2}=0$, then

\begin{eqnarray*}
\varlimsup_{n\rightarrow\infty}\left[\frac{nF(h_n)}{\ln\ln~ n}\right]^{1/2}
N_2=0.
\end{eqnarray*}
}

\subsubsection{Proof of Lemma \ref{lem1}}
Let us set \begin{eqnarray*}
W_{n,i}=\dfrac{1}{F(h_i)^\ell}K\left( \dfrac{\|\chi-{\cal X}_i\|}{h_i}\right)\left[Y_i\mathds{1}_{\left\{|Y_i|\leq b_n\right\}}-r(\chi)\right] & \text{and} & Z_{n,i}= W_{n,i}-\text{E}W_{n,i},
\end{eqnarray*}
and define
\begin{eqnarray*}
S_n=\sum_{i=1}^nZ_{n,i} \ \ \mbox{and} \ \ V_n=\sum_{i=1}^n\text{E}Z_{n,i}^{2}.
\end{eqnarray*}
Observe that
\begin{eqnarray}\label{decomposition2}
V_n&=&\sum\limits_{i=1}^{n}F(h_i)^{-2\ell}\Bigg\{ \text{E}\left(K^2\left(\frac{\|\chi-{\cal X}\|}{h_i} \right)\left[ Y-r(\chi)\right]^2\right)\nonumber\\ \mbox{} & &
+\text{E}\left(K^2\left(\frac{\|\chi-{\cal X}\|}{h_i} \right)Y\left[ 2r(\chi)-Y\right]\mathds{1}_{\left\lbrace \mid Y\mid>b_n \right\rbrace}\right)\Bigg\}\nonumber\\ \mbox{} & &
-\sum\limits_{i=1}^{n}F(h_i)^{-2\ell}\text{E}^2\left(K\left(\frac{\|\chi-{\cal X}\|}{h_i} \right)\left[ Y\mathds{1}_{\left\{|Y|\leq b_n\right\}}-r(\chi)\right]\right)\nonumber\\&:=&A_1+A_2-A_3.
\end{eqnarray}
$A_1$ writes 
\begin{eqnarray*}
A_1&=&\sum\limits_{i=1}^{n}F(h_i)^{-2\ell} \text{E}\left\{K^2\left(\frac{\|\chi-{\cal X}\|}{h_i} \right)\text{E}\left[ (Y-r(\chi))^2|{\cal X}\right]\right\}\\
&=&\sum\limits_{i=1}^{n}\frac{\sigma^2_\varepsilon(\chi) \text{E}K^2\left(\frac{\|\chi-{\cal X}\|}{h_i} \right)}{F(h_i)^{2\ell}} +\frac{\text{E}\left[K^2\left(\frac{\|\chi-{\cal X}\|}{h_i} \right)\{\sigma^2_\varepsilon({\cal X})-\sigma^2_\varepsilon(\chi)\} \right]}{F(h_i)^{2\ell}}\nonumber\\&:=&A_{11}+A_{12}.\label{egalite2}
\end{eqnarray*}
From ${\bf H2}$, by applying Fubini's Theorem, we  have
\begin{eqnarray*}A_{11}=\sum\limits_{i=1}^{n}F(h_i)^{1-2\ell}\sigma^2_\varepsilon(\chi)\left[K^2(1)-\int_0^1(K^2(s))'\tau_{h_i}(s)ds
\right].
\end{eqnarray*}
Applying again Toeplitz's Lemma, ${\bf H3}$ and ${\bf H4}$ allow to get 
 \begin{equation}\label{A11}\frac{A_{11}}{nF(h_n)^{1-2\ell}}\rightarrow\beta_{[1-2\ell]}\sigma^2_\varepsilon(\chi)M_2, \ \ \text{as} \ \  n\rightarrow+\infty.
 \end{equation}
 The second term of  the decomposition of $A_1$ is bounded  using  \eqref{M2}  as 
 \begin{eqnarray*}
A_{12}\leq\sum\limits_{i=1}^{n}F(h_i)^{1-2\ell}\sup\limits_{\chi'\in B(\chi,h_i)}|\sigma^2_\varepsilon(\chi')-\sigma^2_\varepsilon(\chi)|\left[K^2(1)-\int_0^1(K^2(s))'\tau_{h_i}(s)ds
\right], 
\end{eqnarray*}
while the continuity of $\sigma^2_\varepsilon$  (\textbf{H1}) with Toeplitz's Lemma ensure that
\begin{equation}
\label{A12}
\frac{A_{12}}{nF(h_n)^{1-2\ell}}\rightarrow0, \text{ as } n\rightarrow+\infty.
 \end{equation}
Now, let us study the term $A_2$ appearing in the decomposition of $V_n$. Using Cauchy-Schwartz's inequality, and denoting $\left\| K \right\|_{\infty}:=\sup\limits_{t\in [0,1]} K(t)$, we get
\begin{eqnarray*}
A_2&\leq&\|K\|_\infty^2 \sum\limits_{i=1}^{n}F(h_i)^{-2\ell}\left\{\text{E}\left(Y^2\left[2r(\chi)-Y\right]^2\right)P\left(|Y|>b_n\right)\right\}^\frac{1}{2}\\
&\leq&3Q_n\|K\|_\infty^2 \sum\limits_{i=1}^{n}F(h_i)^{-2\ell},
\end{eqnarray*}
where
$Q_n=\left(\max\left\{\Eb\left(Y^4\right),4|r(\chi)|\Eb|Y|^3,4r^2(\chi)\Eb\left(Y^2\right)\right\}P\left(|Y|>b_n\right)\right)^\frac{1}{2}.$\\
We deduce from ${\bf H4}$ and ${\bf H5}$, with the choice $b_n=(\delta \ln n)^{1/\mu}$, that
\begin{equation}\label{A2}
\frac{A_2}{nF(h_n)^{1-2\ell}}=O\left[\frac{e^{-\frac{\lambda b_n^\mu}{2}}\left(\ln n\right)^\frac{2}{\mu}}{F(h_n)}\right]\rightarrow 0, \text{ as }n\rightarrow+\infty \text{ with } \delta>\frac{2}{\lambda}.
\end{equation}
The last term $A_3$ is bounded
\begin{eqnarray*}
|A_3|\leq b_n^2\left[1+o(1)\right]\sum\limits_{i=1}^{n}F(h_i)^{2-2\ell}\left[K(1)-\int_0^1(K'(s))\tau_{h_i}(s)ds\right]^2.
\end{eqnarray*}
It follows from \textbf{H6} that
\begin{equation}\label{A3}
\frac{A_3}{nF(h_n)^{1-2\ell}}=O\left[F(h_n)(\ln n)^\frac{2}{\mu} \right]\rightarrow 0, \text{ as  }n\rightarrow+\infty .
\end{equation}
relations (\ref{A11}), (\ref{A12}), (\ref{A2}) and (\ref{A3}) can be used to derive the following equivalence
\begin{equation}
\label{Vn}
V_n\sim nF(h_n)^{1-2\ell}\beta_{[1-2\ell]}\sigma^2_\varepsilon(\chi)M_2,  \text{  as } n\rightarrow+\infty.
\end{equation}
Next, since $\frac{\ln F(h_n)}{\ln n }\rightarrow0$ as $n\rightarrow+\infty,$  then  the first part of {\bf H6} implies that
\begin{eqnarray*}
\dfrac{nF(h_n)(\ln n)^{-\frac{2}{\mu}}}{\ln\left[nF(h_n)^{1-2\ell}\right]\left\{\ln\ln\left[nF(h_n)^{1-2\ell}\right]\right\}^{2(\alpha+1)}}
\rightarrow\infty, \ \ \text{ as } n \to +\infty.
\end{eqnarray*}
Setting $b_n=(\delta\ln n)^\frac{1}{\mu}$, there  exists $n_0\geq 1$ such that for any  $i\geq n_0$, we have
\begin{eqnarray*}
\dfrac{iF(h_i)(\ln i)^{-\frac{2}{\mu}}}{\ln\left[iF(h_i)^{1-2\ell}\right]\left\{\ln\ln\left[iF(h_i)^{1-2\ell}\right]\right\}^{2(\alpha+1)}}
 >\frac{2\left\| K\right\|_\infty^2\max\left\{|r(\chi)|^2, (\delta\ln i)^\frac{2}{\mu}\right\}}{F(h_i)^{2\ell}}\geq Z_{n,i}^{2}.
 \end{eqnarray*}
Hence, the event ${\left\{Z_{n,i}^2>\frac{iF(h_i)^{1-2\ell}}{\ln \left[iF(h_i)^{1-2\ell}\right]\left\{\ln\ln\left[iF(h_i)^{1-2\ell}\right]\right\}^{2(\alpha+1)}}\right\}}$ is empty for $i\geq n_0.$
 We deduce from (\ref{Vn}) that
 \begin{eqnarray*}
 \sum\limits_{i=1}^\infty\frac{\left(\ln\ln V_i\right)^\alpha}{V_i}\text{E}
 \left(Z_{n,i}^2\mathds{1}_{\left\{Z_{n,i}^2>\frac{V_i}{\ln V_i\left(\ln\ln V_i\right)^{2(\alpha+1)}}\right\}}\right)<\infty.
 \end{eqnarray*}
 Let $S$ be a random function defined on $\left[0,+\infty\right[$ such that for any $t\in[V_n,V_{n+1}[,  S(t)=S_n.$  Using Theorem 3.1 in \cite{Jain}, there exists a Brownian motion $\xi$ such that $$\label{condition4}\left|\frac{S(t)-\xi(t)}{\left(2t\ln\ln t\right)^\frac{1}{2}}\right|= o\left[\left(\ln\ln t\right)^{-\frac{\alpha}{2}}\right] a.s., \ \ \text{as } t \to \infty, \text { for any } t\in[V_n,V_{n+1}[.$$ It follows that
\begin{eqnarray*}
\varlimsup\limits_{t\rightarrow\infty}\frac{S(t)}{\left(2t\ln\ln t\right)^\frac{1}{2}}
=\overline{\lim\limits_{t\rightarrow\infty}}\left[\frac{S(t)-\xi(t)}{\left(2t\ln\ln t\right)^\frac{1}{2}}+\frac{\xi(t)}{\left(2t\ln\ln t\right)^\frac{1}{2}}\right]=1 \ \ a.s.
\end{eqnarray*}
Then we have
\begin{equation}\label{CV1} \frac{S_n}{\sqrt{2V_n\ln\ln V_n}}\rightarrow 1 \ \ a.s., \text{ as $n\rightarrow\infty$},
\end{equation}  by virtue of the definition of $S$  and the fact that  $\frac{V_{n+1}}{V_n}\rightarrow 1$  as $n\to \infty$.
From \eqref{Vn}, we have
\begin{eqnarray*}
\lim\limits_{n\rightarrow\infty}\dfrac{\left\{nF(h_n)^{1-2\ell}\ln\ln\left[nF(h_n)^{1-2\ell}\right]\right\}^{1/2}B_{n,1-\ell}}{
\left(2V_n\ln\ln V_n\right)^\frac{1}{2}}&=&
\frac{\beta_{[1-\ell]}}{\left\{2\beta_{[1-2\ell]}\sigma^2_\varepsilon(\chi)M_2\right\}^{1/2}}.
\end{eqnarray*}
Lemma \ref{lem1} follows from the last convergence and the fact that $S_n=N_1\sum\limits_{i=1}^nF(h_i)^{1-\ell}$, with the help of (\ref{CV1}).\hfill $\square$
\subsubsection{Proof of Lemma \ref{lem3}}
We have
\begin{eqnarray}\label{Bias_as}
N_2
&=&\frac{1}{\sum\limits_{i=1}^nF(h_i)^{1-\ell}}\sum\limits_{i=1}^nF(h_i)^{-\ell}\Bigg\{\Eb\left[ K\left(\frac{\|\chi-{\cal X}\|}{h_i} \right)(r({\cal X})-r(\chi))\right]\nonumber\\\mbox{}&&-\Eb \left[K\left(\frac{\|\chi-{\cal X}\|}{h_i} \right)Y{\mathds{1}}_{\{|Y|>b_n\}}\right]\Bigg\}:= A+B.
\end{eqnarray}
As in the proof of Lemma \ref{Biais1}, we can write
\begin{eqnarray}\label{A_as}
A&=& h_n\frac{\alpha_{[\ell]}}{\beta_{[1-\ell]}}\varphi'(0)M_0\left[1+o(1)\right],
\end{eqnarray}
and then,
\begin{eqnarray*}
\left[\frac{nF(h_n)}{\ln\ln n}\right]^{1/2}A&=&\left[\frac{nF(h_n)}{\ln\ln n}\right]^{1/2}h_n\frac{\alpha_{[\ell]}}{\beta_{[1-\ell]}}\varphi'(0)M_0\left[1+o(1)\right]=o(1),
\end{eqnarray*}
where the last equality follows from the condition $nh_n^2\to 0$. For the second term of the right-hand-side in \eqref{Bias_as},
using Cauchy-Schwartz's inequality and the boundedness of the kernel $K$, we get
\begin{eqnarray*}
|B| &\leq& \frac{\|K\|_{\infty}}{\sum\limits_{i=1}^nF(h_i)^{1-\ell}}\sum_{i=1}^nF(h_i)^{-\ell}\left\{\Eb\left[Y_i^2\right]\Pb\left[|Y_i|>b_n\right]\right\}^{1/2}.
\end{eqnarray*}
From Markov's inequality combined with \eqref{moment-expo}, it follows that
\begin{eqnarray}\label{B_as}
|B| &\leq& \frac{\|K\|_{\infty}}{\sum\limits_{i=1}^nF(h_i)^{1-\ell}}\sum_{i=1}^nF(h_i)^{-\ell}\left\{\Eb\left[Y_i^2\right]\Eb\left[e^{\lambda|Y_i|^\mu}\right]e^{-\lambda b_n^\mu}\right\}^{1/2}\nonumber\\ & = & O\left(\frac{1}{nF(h_n)}\frac{B_{n,-\ell}}{B_{n,1-\ell}}n^{1-\lambda\delta}\left(\ln n\right)^{2/\mu}\right),
\end{eqnarray}
which gives
\begin{eqnarray*}
\left[\frac{nF(h_n)}{\ln\ln n}\right]^{1/2}B&=& O\left(\frac{1}{\sqrt{\ln\ln n}}\frac{1}{\sqrt{nF(h_n)}}n^{1-\lambda\delta}\left(\ln n\right)^{2/\mu}\right)=o(1) \ \ \text{if}\ \ \delta>\frac{1}{\lambda},
\end{eqnarray*}
and Lemma \ref{lem3} is proved. \hfill$\square$
\subsection{Proof of Theorem \ref{normality}}
Using the decomposition   \eqref{egalite1}, we have to show that
\begin{eqnarray}\label{terme1}
\sqrt{nF(h_n)}\left[\varphi_n^{[\ell]}(\chi)-\tilde\varphi_n^{[\ell]}(\chi)\right]\to 0 \ \mbox{a.s.}
\end{eqnarray}
and
\begin{eqnarray}\label{terme2}
\sqrt{nF(h_n)}\left[\tilde\varphi_n^{[\ell]}(\chi)-r(\chi)f_n^{[\ell]}(\chi)\right]\stackrel{\mathcal{D}}{\rightarrow}\mathcal{N}\left(\frac{cM_0\varphi'(0)\alpha_{[\ell]}}{\beta_{[1-\ell]}},\ \frac{\beta_{[1-2\ell]}M_2\sigma_\varepsilon^2(\chi)}{\beta_{[1-\ell]}^2}\right),
\end{eqnarray}
where $\tilde\varphi_n^{[\ell]}$ is defined in \eqref{phi_nl_tilde} and $c$ is such that $\lim\limits_{n\to \infty}\sqrt{nF(h_n)}h_n=c$, since
$f_n^{[\ell]}(\chi)\stackrel{P}{\rightarrow}M_1.$ This later follows  from  the first parts of Lemmas \ref{Biais2} and \ref{Var}.
For \eqref{terme1}, we follow the same lines of proof of \eqref{conv_residual} substituting  $\left[\dfrac{\ln \ln n}{nF(h_n)}\right]^{1/2}$ by $\dfrac{1}{\sqrt{nF(h_n)}}$ to get the desired result. About \eqref{terme2}, using the decomposition \eqref{decomposition11},   it remains to prove Lemmas \ref{variance_term} and \ref{bias_term} below.
\begin{lem} \label{variance_term}
Assume that Assumptions ${\bf H1}-{\bf H5}$ and ${\bf H7}$ hold. Then
\begin{eqnarray*}
\sqrt{nF(h_n)}N_1\stackrel{\mathcal{D}}{\rightarrow}\mathcal{N}\left(0,\frac{\beta_{[1-2\ell]}}{\beta_{[1-\ell]}^2}\sigma_\varepsilon^2(\chi)M_2\right).
\end{eqnarray*}
\end{lem}
\begin{lem} \label{bias_term}
Assume that Assumptions ${\bf H1}-{\bf H5}$ hold. If there exists $c\geq 0$ such that $\lim\limits_{n \to \infty}h_n\sqrt{nF(h_n)}=c$, then
\begin{eqnarray*}
\lim_{n\to \infty}\sqrt{nF(h_n)}N_2=c\frac{\alpha_{[\ell]}}{\beta_{[1-\ell]}}\varphi'(0)M_0.
\end{eqnarray*}
\end{lem}
\subsubsection{Proof of Lemma \ref{variance_term}}
Setting
 \begin{eqnarray*}
W_{n,i}'=\dfrac{\sqrt{nF(h_n)}}{\sum_{i=1}^nF(h_i)^{1-\ell}}W_{n,i} \ \ \text{and} \ \
Z_{n,i}'= W_{n,i}'-\text{E}W_{n,i}',
\end{eqnarray*}
where $W_{n,i}$ is defined in the proof of Theorem \ref{Cv_ps_rnl}, then
\begin{eqnarray*}
\sqrt{nF(h_n)}N_1=\sum_{i=1}^nZ_{n,i}'.
\end{eqnarray*}

To prove Lemma \ref{variance_term}, we first prove that
\begin{eqnarray}\label{condition1}
\lim_{n \to \infty} \sum_{i=1}^n \text{E}(Z_{n,i}^{'2})=\frac{\beta_{[1-2\ell]}}{\beta_{[1-\ell]}^2}\sigma_\varepsilon^2(\chi)M_2,
\end{eqnarray}
and then check that $W_{n,i}'$ satisfies the Lyapounov's condition.
Next, from \eqref{Vn} we have
\begin{eqnarray*}
\sum_{i=1}^n \text{E}(Z_{n,i}^{'2})
&=& \dfrac{nF(h_n)}{\left(\sum_{i=1}^nF(h_i)^{1-\ell}\right)^2}V_n=
\frac{1}{nF(h_n)^{1-2\ell}}\dfrac{1}{B_{n,1-\ell}^2}V_n\\&=& \frac{\beta_{[1-2\ell]}}{\beta_{[1-\ell]}^2}\sigma_\varepsilon^2(\chi)M_2\left[1+o(1)\right]
\end{eqnarray*}
which proves \eqref{condition1}.
To check the Lyapounov's condition, set $p>2$, we have
\begin{eqnarray*}
\sum_{i=1}^n\Eb\left(|Z_{n,i}'|^p\right)=\sum_{i=1}^n\Eb\left(|Z_{n,i}'|^{p-2}Z_{n,i}'^2\right).
\end{eqnarray*}
Since
\begin{eqnarray*}
\left|W_{n,i}'\right|&\leq & \dfrac{\|K\|_{\infty}\sqrt{nF(h_n)}}{\sum_{i=1}^nF(h_i)^{1-\ell}}F(h_i)^{-\ell}|b_n-r(\chi)|,
\end{eqnarray*}
it follows that
\begin{eqnarray*}
\sum_{i=1}^n\Eb\left(|Z_{n,i}'|^p\right)
\leq \dfrac{(nF(h_n))^{\frac{p}{2}}\sum\limits_{i=1}^nF(h_i)^{-p\ell}Var\left(K\left(\frac{\|\chi-\X_i\|}{h_i}\right)\left(Y_i\mathds{1}_{\left\{|Y_i|\leq b_n\right\}}-r(\chi)\right)\right)}{2^{2-p}\|K\|_{\infty}^{2-p}\left|b_n-r(\chi)\right|^{2-p}\left(\sum\limits_{i=1}^nF(h_i)^{1-\ell}\right)^{p}}.
\end{eqnarray*}
Using the same decomposition as in \eqref{decomposition2}, we have

\begin{eqnarray}\label{repet}
\lefteqn{\sum_{i=1}^nF(h_i)^{-p\ell}Var\left(K\left(\frac{\|\chi-\X_i\|}{h_i}\right)\left(Y_i\mathds{1}_{\left\{|Y_i|\leq b_n\right\}}-r(\chi)\right)\right)}\nonumber\\&=&
\sum\limits_{i=1}^{n}F(h_i)^{-p\ell}\Bigg\{ \Eb\left( K^2\left(\frac{\|\chi-{\cal X}\|}{h_i} \right)\left[ Y-r(\chi)\right]^2\right)\nonumber\\ \mbox{}&&
+\text{E}\left(K^2\left(\frac{\|\chi-{\cal X}\|}{h_i} \right)Y\left[ 2r(\chi)-Y\right]\mathds{1}_{\left\lbrace \mid Y\mid>b_n \right\rbrace}\right)\Bigg\}\nonumber\\ \mbox{} & &
-\sum\limits_{i=1}^{n}\frac{\text{E}^2K\left(\frac{\|\chi-{\cal X}\|}{h_i} \right)\left[ Y\mathds{1}_{\left\{|Y|\leq b_n\right\}}-r(\chi)\right]}{F(h_i)^{p\ell}}:=B_1+B_2-B_3.
\end{eqnarray}
Setting $b_n=(\delta\ln n)^{\frac{1}{\mu}}$ for some $\delta,\mu>0$ and following the same lines as in the proof of \eqref{A11}, \eqref{A12}, \eqref{A2} and \eqref{A3} substituting the exponent $2$ by $p$ in all the expressions, we have
\begin{eqnarray*}
B_1&=&O\left(nF(h_n)^{1-p\ell}\right).
\end{eqnarray*}
From Toeplitz's Lemma, we can write
\begin{eqnarray*}
\dfrac{(nF(h_n))^{\frac{p}{2}}}{\left(\sum_{i=1}^nF(h_i)^{1-\ell}\right)^{p}}\left|b_n-r(\chi)\right|^{p-2}B_1&=& O\left(\left(\frac{(\ln n)^{\frac{1}{\mu}}}{\sqrt{nF(h_n)}}\right)^{p-2}\right) = o(1).
\end{eqnarray*}
Next, for the second expression $B_2$ of \eqref{repet}, we get
\begin{eqnarray*}
\frac{B_2}{nF(h_n)^{1-p\ell}}&=&O\left(\frac{\exp\left(-\frac{\lambda b_n^\mu}{2}\right)\left(\ln n\right)^{\frac{2}{\mu}}}{F(h_n)}\right)=o(1) \ \text{with} \ \ \delta>\frac{2}{\lambda}.
\end{eqnarray*}
It follows again from Toeplitz's Lemma that
\begin{eqnarray*}
\dfrac{(nF(h_n))^{\frac{p}{2}}}{\left(\sum_{i=1}^nF(h_i)^{1-\ell}\right)^{p}}\left|b_n-r(\chi)\right|^{p-2}B_2&=& O\left(\left(\frac{(\ln n)^{\frac{1}{\mu}}}{\sqrt{nF(h_n)}}\right)^{p-2}\right) = o(1).
\end{eqnarray*}
In the same manner as in the proof of \eqref{A3}, we have
\begin{eqnarray*}
\frac{B_3}{nF(h_n)^{1-p\ell}}=O\left[F(h_n)(\ln n)^\frac{2}{\mu} \right].
\end{eqnarray*}
Then
\begin{eqnarray*}
\dfrac{(nF(h_n))^{\frac{p}{2}}}{\left(\sum_{i=1}^nF(h_i)^{1-\ell}\right)^{p}}\left|b_n-r(\chi)\right|^{p-2}B_3&=&O\left(\left(\frac{(\ln n)^{\frac{1}{\mu}}}{\sqrt{nF(h_n)}}\right)^{p-2}\right) = o(1),
\end{eqnarray*}
which concludes the proof of Lemma \ref{variance_term}.\hfill $\square$

\subsubsection{Proof of Lemma \ref{bias_term}}
We go back to the decomposition of \eqref{Bias_as} in the proof of Lemma \ref{lem3}.\\
On one hand, from \eqref{A_as}, we write
\begin{eqnarray}
\sqrt{nF(h_n)}A=\sqrt{nF(h_n)}h_n\frac{\alpha_{[\ell]}}{\beta_{[1-\ell]}}\varphi'(0)M_0\left[1+o(1)\right]\label{bias_term_nor1}.
\end{eqnarray}
On the other hand from \eqref{B_as}, we get
\begin{eqnarray}\label{bias_term_nor2}
\sqrt{nF(h_n)}B&=&O\left(\frac{1}{\sqrt{nF(h_n)}}\right)=o(1) \ \ \text{if}\ \ \delta>\frac{1}{\lambda},
\end{eqnarray}
and Lemma \ref{bias_term}  follows from the combination of \eqref{bias_term_nor1} and \eqref{bias_term_nor2}. \hfill$\square$
\subsubsection{Proof of Corollary \ref{normality2}}
From the standard Glivenko-Cantelli type results, we have 
\begin{eqnarray}\label{glivenko}
\frac{F(h_n)}{\widehat{F}(h_n)}\stackrel{\mathbb{P}}{\rightarrow}1.
\end{eqnarray}
Then,  to prove Corollary \ref{normality2}, it suffices to show that $\widehat\beta_{[1]}\stackrel{\mathbb{P}}{\rightarrow}\beta_{[1]}$. To see this, we use the following result $$\sup_{x\in\mathbb{R}}\left|\widehat{F}(x)-F(x)\right|=O((\ln n/n)^{1/2}),$$ established by  \cite{D77} and write
\begin{eqnarray*}
\widehat\beta_{[1]}-\beta_{[1]}&=&\frac{1}{n\widehat{F}(h_n)}\sum_{i=1}^n\left(\widehat{F}(h_i)-F(h_i)\right)+\frac{F(h_n)}{\widehat{F}(h_n)}B_{n,1}-\beta_{[1]}\\ &\leq&\frac{F(h_n)}{\widehat{F}(h_n)}\frac{1}{nF(h_n)}\sup_{x\in\mathbb{R}}\left|\widehat{F}(x)-F(x)\right|+\left|\frac{F(h_n)}{\widehat{F}(h_n)}B_{n,1}-\beta_{[1]}\right|=o_{\mathbb P}(1),
\end{eqnarray*}
by virtue of {\bf H4(ii)} and \eqref{glivenko}. \hfill$\square$


\begin{thebibliography}{99}
 \bibitem[Ahmad and Lin(1976)]{AL76} Ahmad, I. and Lin, P.E. (1976). Nonparametric sequential estimation of a multiple regression function. \textit{Bull. Math. Statist.}, \textbf{17}, 63-75.
\bibitem[Amiri(2012)]{Amiri} Amiri, A. (2012). Recursive regression estimators with application to nonparametric prediction. \emph{J. Nonparametr. Stat.}  {\bf 24} (1),   169-186.
 \bibitem[ Besse et al.(1997)]{BCF97} Besse, P., Cardot, H. and Ferraty, F. (1997). Simultaneous nonparametric regressions of unbalanced longitudinal data. \textit{Comp. Statist. Data Anal.}, \textbf{24}, 255-270.
\bibitem[Bosq and Cheze-Payaud(1999).]{bosq-cheze} Bosq, D., and Cheze-Payaud, N. (1999). Optimal asymptotic quadratic error of nonparametric regression function estimates for a continuous-time process from sampled-data. \emph{Statistics}, {\bf 32} (3), 229--247.
\bibitem[Cardot et al.(2003)]{CFS03} Cardot, H., Ferraty, F. and Sarda, P. (2003). Splines estimators for the functional linear model. \textit{Statistica Sinica}, \textbf{13}, 571-591.
\bibitem[Collomb(1976)]{Collomb1976} Collomb, G. (1976). Estimation nonparam\'etrique de la r\'egression par la m\'ethode du noyau.  {\it Unpublished PhD Thesis},  Universit\'e Paul Sabatier, Toulouse, France.
\bibitem[Crambes et al.(2009)]{CKS09} Crambes, C., Kneip, A. and Sarda, P. (2009). Smoothing splines estimators for functional linear regression, \textit{Ann. Statist.}, \textbf{37}, 35-72.
 \bibitem[Devroye (1977)]{D77} Devroye, L., (1977). A uniform bound for the deviation of empirical distribution functions. \textit{J. Mult. Anal.}, \textbf{7}, 594-597.
\bibitem[Devroye and Wagner(1980)]{DW80} Devroye, L. and Wagner, T.J. (1980). Distribution-free consistency results in nonparametric discrimination and regression function estimation. \textit{Ann. Statist.}, \textbf{8}, 231-239.
 \bibitem[Ferraty et al.(2007)]{Ferraty2007} Ferraty, F. Mas, A. and Vieu, P. (2007). Nonparametric regression on functional data: Inference and practical aspects. \emph{Aust. N. Z. J. Stat.,} {\bf 49} (3), 267--286.
 \bibitem[Ferraty and Romain(2010)]{FerratyRomain} Ferraty, F. and Romain, Y. (2010). Handbook on functional data analysis and related fields. {\it  Oxford University Press}, Oxford.
 \bibitem[Ferraty and Vieu(2006)]{Ferraty2006} Ferraty, F. and Vieu, P. (2006). Nonparametric Modelling for Functional Data. Methods, Theory, Applications and Implementations. {\it Springer-Verlag},  London.
 \bibitem[Frank and Friedman (1993)]{FrankFriedman} Frank, I.E. and Friedman, J.H. (1993). A statistical view of some chemometrics regression tools. \textit{Technometrics}, \textbf{35}, 109--135.
 \bibitem[Gonz\`alez Manteiga and Vieu(2007)]{MantVieu} Gonz\`alez Manteiga, W. and Vieu, P. (2007). Statistics for functional data. \textit{Comp. Stat. and Data Anal.}, \textbf{51}, 4788-4792.
 \bibitem[Jain et al.(1975)]{Jain} Jain, C., Jogdeo, K. and Stout, W. (1975). Upper and lower functions for martingales and mixing processes, \emph{Ann. Probab.} { \bf 3}  (1), 119--145.
 \bibitem[Ling and Wu (2012)]{LINGWU} Ling, N.  and   Wu, Y. (2012). Consistency of modified kernel regression estimation for functional data, \emph{Statistics} {\bf 46} (2), 149--158.
 \bibitem[Masry(1986)]{Masry86} Masry, E. (1986). Recursive probability density estimation for weakly dependent
stationary processes, {\it IEEE Trans. Inform. Theory} {\bf 32} (2), 254--267.
 \bibitem[Masry(2005)]{Masry2005} Masry, E. (2005). Nonparametric regression estimation for dependent functional data: Asymptotic normality, \emph{Stochastic Process. Appl.} {\bf 115} 155--177.
 \bibitem[Rachdi and Vieu(2007)]{RachdiVieu} Rachdi, M. and Vieu, P. (2007). Nonparametric regression for functional data: automatic smoothing parameter selection. \textit{J. Statist. Plann. Inference}, \textbf{137}, 2784-2801.
 \bibitem[Ramsay and Dalzell(1991)]{RamsayDalzell} Ramsay, J.O. and Dalzell, C.J. (1991). Some tools for functional data analysis (with discussion), \textit{J. Roy. Statist. Soc., Ser. B}, \textbf{53}, 539-572.
 \bibitem[Ramsay  and Silverman(2002)]{RamsaySilverman02} Ramsay, J.O. and Silverman, B.W. (2002). Applied functional data analysis. {\it Springer-Verlag}, New-York.
\bibitem[Ramsay and Silverman(2005)]{RamsaySilverman05} Ramsay, J.O. and Silverman, B.W. (2005). \textit{Functional data analysis ($2^{\text{nd}}$ Ed.)}. Springer-Verlag, New-York.
  \bibitem[Roussas(1992)]{Roussas92} Roussas, G. G.(1992). Exact rates of almost sure convergence of a recursive kernel estimate of a probability densiy function: Application to regression and hazard rate estimation, \emph{J. Nonparametr. Stat.,} {\bf 1} (3), 171-- 195.
 \end{thebibliography}
\end{document}